\documentclass{amsart}    
\usepackage{amsmath} 
\usepackage{paralist}
\usepackage{graphics} 
\usepackage{epsfig} 
\usepackage{graphicx}  \usepackage{epstopdf} 
\usepackage[colorlinks=true]{hyperref}
\hypersetup{urlcolor=blue, citecolor=red}

\newtheorem{Proposition}{Proposition}[section]

\newtheorem{Lemme}{Lemma}[section]
\newtheorem{Theoreme}{Theorem}[section]

\newtheorem{Corollaire}{Corollary}[section]
\newtheorem{Remarque}{Remark}

\def \f{\vec{f}} 

\def \u{\vec{u}}
\def \v{\vec{v}}

\def \P{\mathbb{P}}

\def \R{\mathbb{R}}

\def \Rt{\mathbb{R}^3}
\def \finpv{\hfill $\blacksquare$}
\def \pv{{\bf{Proof.}}~} 

\def \ds{\displaystyle}

\usepackage{amssymb}

\title[\bf An $L^p$-theory for fractional stationary Navier-Stokes equations] %
{ An $L^p$-theory for fractional stationary Navier-Stokes equations} 
\author[ Oscar Jarr\'in and Gast\'on Vergara-Hermosilla]{}


\email{oscar.jarrin@udla.edu.ec}
\email{gaston.vergarahermosilla@univ-evry.fr} 

\thanks{$^*$Corresponding author:  Oscar Jarr\'in}

\begin{document}
	\maketitle
	
	\centerline{\scshape Oscar Jarr\'in$^*$}
	\medskip
	{\footnotesize
		\centerline{Escuela de Ciencias Físicas y Matemáticas}
		\centerline{Universidad de Las Américas}
		\centerline{V\'ia a Nay\'on, C.P.170124, Quito, Ecuador}
	} 
	
	\medskip
	
	\centerline{\scshape Gast\'on Vergara-Hermosilla}
	\medskip
	{\footnotesize
		\centerline{Laboratoire de Mathématiques et Modélisation d'Évry, CNRS UMR 8071}
		\centerline{Université Paris-Saclay} 
		\centerline{91025, \'Evry, France} 
	}

	\bigskip

\begin{abstract}
We consider  the  stationary (time-independent) Navier-Stokes equations in the whole three-dimensional space, under  the action of a source term and with the fractional Laplacian operator $(-\Delta)^{\alpha/2}$ in the diffusion term. In the framework of Lebesgue and Lorentz spaces, we find some natural sufficient conditions on the external force and on the parameter $\alpha$   to prove the existence and in some cases nonexistence  of solutions. Secondly, we obtain sharp pointwise decay rates   and asymptotic profiles of solutions, which strongly depend on $\alpha$. Finally,  we also prove the global regularity of solutions. As a bi-product, we obtain some uniqueness theorems so-called Liouville-type results. On the other hand, our regularity result yields a new regularity criterion for the classical ({\it i.e. }with $\alpha=2$) stationary Navier-Stokes equations.   
\end{abstract}
{\footnotesize \tableofcontents}
\section{Introduction}
This paper  considers  the 3D incompressible, stationary  and fractional  Navier-Stokes equations in the whole space $\Rt$:
\begin{equation*}
(-\Delta)^{\frac{\alpha}{2}} \u + (\u \cdot \vec{\nabla}) \u + \vec{\nabla}P= \f, \qquad  \text{div}(\u)=0, \quad \alpha>0.
\end{equation*}
Here, the functions $\u:\Rt \to \Rt$ and $P:\Rt\to \R$ are the velocity and the pressure of the fluid respectively, while $\f:\Rt \to \Rt$ represents a given external force acting on this equation. Moreover, with a minor loss of generality, we set the viscosity constant equal to one. 

\medskip

The divergence-free property of $\u$ yields to easily deduce the pressure $P$ from the velocity $\u$ and the external force $\f$ by the expression
\begin{equation}\label{Pressure}
P= \frac{1}{-\Delta} \text{div}\big( (\u \cdot \vec{\nabla}) \u - \f \big).
\end{equation}
This fact allows us to focus our study on the velocity $\u$ and to consider the equations:
\begin{equation}\label{Frac-NS}
(-\Delta)^{\frac{\alpha}{2}} \u + \P \big( (\u \cdot \vec{\nabla}) \u \big) = \P(\f), \qquad  \text{div}(\u)=0, \quad \alpha>0,
\end{equation}
where $\P$ stands for the Leray's projector. 

\medskip

One of the main features of   equation (\ref{Frac-NS})  is the fractional Laplacian operator $(-\Delta)^{\frac{\alpha}{2}}$  in the diffusion term. Recall that this operator is defined at  the Fourier level by the symbol $|\xi|^\alpha$, whereas in the spatial variable we have 
\[ (-\Delta)^{\alpha/2} \u(t,x)= C_{\alpha}\, {\bf p.v.} \int_{\Rt} \frac{\u(t,x)- \u(t,y)}{|x-y|^{3+\alpha}} dy,  \]
where $C_\alpha>0$ is a constant depending on $\alpha$, and  ${\bf p.v.}$ denotes the principal value.

\medskip

Experimentally, this operator  has been successfully employed to model anomalous reaction-diffusion process
in porous media models \cite{meerschaert1999pseudo,meerschaert2001pseudo} and in computational turbulence models \cite[Chapter 13.2]{pope2003pseudo}. In these
last models, the term $(-\Delta)^{\frac{\alpha}{2}}$   is used to characterize anomalous viscous diffusion effects in turbulent fluids which are driven by the parameter $\alpha$.

\medskip

Mathematically, in the setting of bounded and smooth domains $\Omega\subset \Rt$, linear and nonlinear  fractional elliptic equations have been extensively studied in the setting of $L^p$-spaces, see for instance \cite{abramyan1994pseudo,bhatka2019pseudo,caffarelli2016pseudo,chamorropoggi2023p,dipierro2017pseudo,dong2012pseudo,muratori2016p,tang2016new} and the references therein. Motivated by these works, this article is devoted to develop a theory for the equation (\ref{Frac-NS}) in the framework of  Lebesgue spaces and in the more general setting of  Lorentz spaces.  Our main objective is to understand the effects of the  parameter $\alpha$ in the qualitative study of this equation, precisely, the existence and nonexistence of $L^p$-solutions and $L^{p,q}$-solutions, pointwise decaying properties,   asymptotic profiles  and  regularity properties. 
This seems to be, to the best of our knowledge, the first application of these functional spaces in the analysis of fractional stationary  Navier-Stokes equations in the whole three-dimensional space. 
Moreover,  it is worth highlighting the previously cited works are not longer
valid in $\Rt$ due to the lack of some of their key tools, for instance, the embedding properties of  $L^p(\Omega)$-
spaces and  compact Sobolev embeddings. We thus use different approaches to study these qualitative properties. 

\medskip 

In the case $\alpha=2$, equations (\ref{Frac-NS}) coincide with  the classical stationary Navier-Stokes equations
\begin{equation}\label{Classical-NS}
-\Delta  \u + \P \big( (\u \cdot \vec{\nabla}) \u \big) =  \P(\f), \quad \text{div}(\u)=0.
\end{equation}
The $L^p$-theory for these equations was studied in \cite{bjorland2010p}, where the authors exploit powerful tools of Lorentz spaces to study some of the aforementioned qualitative properties, mainly the existence and asymptotic profiles of solutions.

\medskip

In this paper, we generalize some of their results to the \emph{fractional case} of equation  
 (\ref{Frac-NS}) and we highlight some interesting \emph{phenomenological effects}, principally due to the fractional Laplacian operator. Going further,  we also develop a new approach to study the global regularity of solutions. This result has two main consequences, on the one hand, for the classical case (when $\alpha=2$) we obtain a new regularity criterion for weak solutions to  the equation (\ref{Classical-NS}) and, on the other hand,  for the fractional case (when $\frac 5 3<\alpha<2$) our regularity result yields some uniqueness properties  of solutions to the  homogeneous (with $\f\equiv 0$) equation   (\ref{Frac-NS}). This last result is of independent interest and it is also known as the \emph{Liouville}-type  problem for the \emph{fractional} Navier-Stokes equations.   

\medskip

\noindent{\bf Main results}. Recall that for a measurable function  $f: \Rt \to \R$  and for a parameter $\lambda \geq 0$ we define  the distribution function
\[ d_f(\lambda)= dx \left( \left\{  x \in \Rt: \ | f(x)|>\lambda \right\} \right), \]
where $dx$ denotes the Lebesgue measure. Then, the re-arrangement function $f^{*}$ is defined by the expression 
\[ f^{*}(t)= \inf \{ \lambda \geq 0: \ d_f(\lambda) \leq t\}. \]
By definition, for $1\leq p <+\infty$ and $1\leq q \leq +\infty$ the Lorentz space $L^{p,q}(\Rt)$ is the space  of measure functions $f: \Rt \to \R$ such that $\| f \|_{L^{p,q}}<+\infty$, where:
\begin{equation*}
\| f \|_{L^{p,q}}=
\begin{cases}\vspace{3mm}
\ds{\frac{q}{p} \left( \int_{0}^{+\infty} (t^{1/p} f^{*}(t))^{q} dt\right)^{1/q}}, \quad q<+\infty, \\
\ds{\sup_{t>0} \, t^{1/p} f^{*}(t)}, \quad q=+\infty.
\end{cases}
\end{equation*}
It is worth mentioning some important properties of these spaces. The quantity $\| f \|_{L^{p,q}}$ is often used as a norm, even thought it does not verify the triangle inequality. However, there exists  an equivalent norm (strictly speaking) which makes these spaces into Banach spaces. On the other hand, these spaces are homogeneous of degree $-\frac{3}{p}$ and for $1 \leq q_1 <  p < q_2\leq +\infty$ we have the continuous embedding 
\[ L^{p,q_1}(\Rt) \subset L^{p}(\Rt)=L^{p,p}(\Rt) \subset L^{p,q_2}(\Rt). \]  
Finally, for $p=+\infty$ we  have the identity  $L^{\infty,\infty}(\Rt)=L^\infty(\Rt)$.

\medskip

In our first result, we find some conditions on the external force $\f$ and a range of values of the parameter $\alpha$  to construct $L^p$-weak  solutions and $L^{p,q}$-weak solutions to equation (\ref{Frac-NS}). For this, we shall perform the following program: at point $(A)$ below, first  we solve this equation in the Lorentz space $L^{\frac{3}{\alpha-1},\infty}(\Rt)$. This particular space  is invariant under the natural scaling of equations (\ref{Frac-NS}): $(\u,P,\f)\mapsto (\u_\lambda, P_\lambda, \f_\lambda)$, where
\[\u_\lambda(x)=\lambda^{\alpha-1}\u(\lambda x), \quad P_\lambda(x)=\lambda^{2\alpha-2}P(\lambda x) \quad \mbox{and} \quad  \f_\lambda(x)=\lambda^{2\alpha-1}\f(\lambda x), \quad \lambda>0,\]
and this fact   allows us to  apply the Picard iterative scheme to the fixed point problem:
\begin{equation}\label{Frac-NS-fixed-point}
 \u =  - (-\Delta)^{-\frac{\alpha}{2}}  \P (\text{div}(\u \otimes \u)) + (-\Delta)^{-\frac{\alpha}{2}} \P (\f).
\end{equation}
Precisely, by a sharp study of the kernel associated to the operator $(-\Delta)^{-\frac{\alpha}{2}}  \P (\text{div}(\cdot))$ 
we are able to prove the estimate $\| (-\Delta)^{-\frac{\alpha}{2}}  \P (\text{div}(\u \otimes \u)) \|_{L^{\frac{3}{\alpha-1},\infty}} \lesssim \| \u \|^{2}_{L^{\frac{3}{\alpha-1}, \infty}}$, provided that $1<\alpha<\frac{5}{2}$. 

\medskip

To the best of our knowledge, this well established procedure does not work for $0<\alpha<1$ nor $\frac{5}{2}\leq \alpha$, and these cases will be matter of (far from obvious) future research. Consequently, from now on, we set $1<\alpha < \frac{5}{2}$. 

\medskip

Continuing with our program, at points $(B)$ and $(C)$ below we prove some \emph{persistence} properties of $L^{\frac{3}{\alpha-1}}$-solutions. Specifically, we find a range of values for the parameter $p$ for which these solutions also belong to the $L^p$-space or the $L^{p,q}$-space, as long as the external force verifies additional suitable conditions.    
\begin{Theoreme}\label{Th1} Let $1<\alpha < \frac{5}{2}$. Assume that  $(-\Delta)^{-\frac{\alpha}{2}} \f \in L^{\frac{3}{\alpha-1},\infty}(\Rt)$. There exists a universal quantity $\varepsilon_0(\alpha)>0$, which only depends on $\alpha$, such that if 
\[ \| (-\Delta)^{-\frac{\alpha}{2}} \P ( \f ) \|_{L^{\frac{3}{\alpha-1},\infty}} < \varepsilon_0(\alpha),\] 
then the following statements hold:
\begin{enumerate}
    \item[(A)] The  equation (\ref{Frac-NS}) has a solution $\u \in L^{\frac{3}{\alpha-1},\infty}(\Rt)$ satisfying and uniquely defined by the condition $\| \u \|_{L^{\frac{3}{\alpha-1},\infty}} \leq 2 \| (-\Delta)^{-\frac{\alpha}{2}} \P ( \f ) \|_{L^{\frac{3}{\alpha-1},\infty}}$. 

    \medskip 
 
    \item[(B)] Let $\frac{3}{4-\alpha}<p<+\infty$. Assume that $(-\Delta)^{-\frac{\alpha}{2}} \f \in  L^p(\Rt)$,  then the solution constructed above  verifies $\u \in L^p(\Rt)$.  Moreover, this fact holds in the space $L^{p,q}(\Rt)$ with  $1\leq q \leq +\infty$.

    \medskip
    
    \item[(C)] Finally, for the end point $p=\frac{3}{4-\alpha}$, assume that   $(-\Delta)^{-\frac{\alpha}{2}} \f  \in L^{\frac{3}{4-\alpha}, \infty}(\Rt)$. Then  the solution constructed at point  $(A)$ verifies  $\u \in L^{\frac{3}{4-\alpha}, \infty}(\Rt)$. Moreover, for the end point $p=+\infty$ it holds $\u \in L^\infty(\Rt)$, provided that $(-\Delta)^{-\frac{\alpha}{2}} \P (\f) \in L^{\infty}(\Rt)$.  
\end{enumerate}
\end{Theoreme}

Some comments are in order. As mentioned, at point $(B)$ we show the existence of $L^p$-solutions and $L^{p,q}$-solutions (with $1\leq q \leq +\infty$) for the range of values $\frac{3}{4-\alpha}<p<+\infty$, while, at point $(C)$,  we prove the existence of $L^{\frac{3}{4-\alpha},\infty}$-solutions and $L^\infty$-solutions. The main remark is that these results are \emph{optimal}, and later in Theorem \ref{Prop-non-existence} below we prove the \emph{non-existence} of $L^{p,q}$-solutions for   $1\leq p \leq \frac{3}{4-\alpha}$ and  $1\leq q <+\infty$.

\medskip
In the setting of Sobolev spaces,  existence of  weak  solutions to equation (\ref{Frac-NS}) has been proven in \cite{chamorropoggi2023p} and \cite{tang2016new}.  In the first work, for a divergence-free external force $\f \in \dot{H}^{-1}(\Rt)\cap \dot{H}^{-\frac{\alpha}{2}}(\Rt)$, the authors use the Schaefer fixed point theorem to obtain a weak $\dot{H}^{\frac{\alpha}{2}}$-solution as the limit of a regularized problem. This approach allow them to consider the range $0<\alpha < 2$, but they lose uniqueness. On the other hand, in the second work \cite{tang2016new},  the authors perform a Galerkin-type method to construct a weak suitable $\dot{H}^{\frac{\alpha}{2}}$-solution with $1<\alpha<2$. As before, their uniqueness is still an open problem.
Note that, the existence of a unique stationary solution for the forced tridimensional fractional Navier-Stokes-Coriolis equations,  has been recently proved in \cite{MR4587606} by considering Fourier-Besov spaces. This result remains for any arbitrarily large external force, whenever the Coriolis parameter is taken large enough.

\medskip

Our approach is completely different  and it strongly exploits  the structure of equation (\ref{Frac-NS-fixed-point}).  This approach seems to have some advantages in the study of equation (\ref{Frac-NS}). First, as was already mentioned, we are able to obtain optimal results when studying the existence of $L^p$-solutions. On the other hand, the main assumption of the external force: $(-\Delta)^{-\frac{\alpha}{2}}\P (\f) \in L^{\frac{3}{\alpha-1},\infty}(\Rt)$ allows us to consider a variety of  external forces, for instance, when $\alpha=2$, one can considers a very rough external force $\f=  (c_1\delta_0, c_2\delta_0, c_3\delta_0)$, where $\delta_0$ denotes the Dirac mass at the origin and $c_1,c_2,c_3\in \R$ are suitable numerical constants.
Moreover, when $\alpha\neq 2$, one can consider homogeneous external forces   $| \f (x) | \sim \frac{c}{|x|^{2\alpha-4}}$,  with $x\neq 0$ and $c>0$. Finally, in contrast to \cite{chamorropoggi2023p,tang2016new}, equation (\ref{Frac-NS-fixed-point}) allows us to study sharp asymptotic profiles and pointwise decay rates  of solutions to  equation (\ref{Frac-NS}). 

\medskip

In our next result, for $\theta \geq 0$,   we use the homogeneous weighted $L^\infty$-space
\begin{equation*}
L^{\infty}_{\theta}(\Rt)= \big\{ f \in L^{\infty}_{loc}(\Rt \setminus \{0\}): \ \| f \|_{L^{\infty}_{\theta}} = \mbox{ess sup}_{x \in \Rt} | x |^\theta |f(x)|<+\infty \big\}, 
\end{equation*}
where we have the identity $L^{\infty}_{0}(\Rt)=L^{\infty}(\Rt)$. 
Similar weighted $L^\infty$-spaces has been considered in \cite{MR4518052} to study the solvability and asymptotic profile of solutions for the instationary fractional incompressible Navier–Stokes equations in $\Rt$, however, to the best of the authors knowledge, they are no been used in the analysis of the steady case of these equations. 
Motivated by this, and  by recalling the continuous embedding $L^{\infty}_{\alpha-1}(\Rt) \subset L^{\frac{3}{\alpha-1},\infty}(\Rt)$, we note that it is natural to solve equation (\ref{Frac-NS-fixed-point}) in the smaller \emph{scale invariant} space $L^{\infty}_{\alpha-1}(\Rt)$, for small external forces verifying $(-\Delta)^{-\frac{\alpha}{2}}\P (\f) \in L^{\infty}_{\alpha-1}(\Rt)$. Thereafter, by following the ideas of Theorem \ref{Th1}, we find a range of values for the parameter $\theta$ for which pointwise decaying properties on the external force are propagated to solutions.  
\begin{Theoreme}\label{Th2} Let $1<\alpha<\frac{5}{2}$. Assume that the external force $\f$ verifies $(-\Delta)^{\frac{\alpha}{2}} \P (\f) \in L^{\infty}_{\alpha-1}(\Rt)$. There exists a universal quantity $0<\varepsilon_1(\alpha)<\varepsilon_0(\alpha)$, which only depends on $\alpha$, such that if
\begin{equation*}
\| (-\Delta)^{\frac{\alpha}{2}} \P (\f)  \|_{L^{\infty}_{\alpha-1}} < \varepsilon_1(\alpha),
\end{equation*}
then the following statements hold:

\medskip

\begin{enumerate}
    \item[(A)]  The  equation (\ref{Frac-NS}) has a solution $\u \in L^{\infty}_{\alpha-1}(\Rt)$ satisfying and uniquely defined by the condition $\| \u \|_{L^{\infty}_{\alpha-1}} \leq 2 \| (-\Delta)^{-\frac{\alpha}{2}} \P ( \f ) \|_{L^{\infty}_{\alpha-1}}$. 
    
      \medskip
 
    \item[(B)] If $0\leq \theta \leq  4 - \alpha$ and 
    $(-\Delta)^{-\frac{\alpha}{2}} \mathbb{P} (\f) \in L^{\infty}_{\theta}(\Rt)$ the solution obtained above holds  $\u \in L^{\infty}_{\theta}(\Rt)$. 
\end{enumerate}
\end{Theoreme}

Pointwise decaying properties of solutions allow us to deduce some interesting facts, which we shall explain  in the following  corollaries and propositions. 

\medskip 

First, by point $(B)$ above we directly obtain  the following estimate, which shows an explicit  decay rate of solutions to equation (\ref{Frac-NS}).
\begin{Corollaire}\label{Corolary-pointwise-estimate} Within the framework of Theorem \ref{Th2}, assume that $(-\Delta)^{\frac{-\alpha}{2}} \mathbb{P} ( \f)  \in L^\infty(\Rt)\cap L^\infty_{4-\alpha}(\Rt)$. Then, there exist a numerical constant $C>0$ such that
\begin{equation}\label{sep_1_gaston}
    |\u(x)| \leq \frac{C}{1 + |x|^{4-\alpha}}.  
\end{equation}
\end{Corollaire}

As noticed, this decay rate is driven by the parameter $\alpha$ in the fractional Laplacian operator in equation (\ref{Frac-NS}). Moreover, the main remark is that this decay rate is also \emph{optimal} and in Theorem \ref{Prop-non-existence} below we show that solutions to equation (\ref{Frac-NS}) cannot decay at infinity faster than $\frac{1}{|x|^{4-\alpha}}$.

\medskip

Both non-existence  results in Lebesgue and Lorentz spaces as well as  the optimality of the  decay rate above are obtained by sharp asymptotic profiles of solutions to equations (\ref{Frac-NS}).  For this, we remark that the kernel associated to the operator $(-\Delta)^{-\frac{\alpha}{2}}\P(\text{div}(\cdot))$ is obtained as a tensor $m_\alpha= (m_{i,j,k})_{1\leq i,j,k\leq 3}$, where $m_{i,j,k}(x)$ is a homogeneous function of order $4-\alpha$ and smooth outside the origin (see Lemma \ref{Lem-Kernel} below for all the details). Moreover, we recall that  for $i=1,2,3$ the term $ \ds{m_{\alpha} (x) : \int_{\Rt} (\u \otimes \u )(y)dy}$ is defined by 
\begin{equation}\label{Def-two-dots}
    \left(
m_{\alpha} (x) : \int_{\Rt} (\u \otimes \u )(y)dy
    \right)_i
    =
    \sum_{j,k=0}^3
m_{i,j,k}(x)
 \left(
\int_{\Rt} u_j (y)  u_k (y)dy
    \right).
\end{equation}

\begin{Proposition}\label{Proposition-asymptotic-profile}
Under the same hypotheses of Theorem \ref{Th2}, assume that $(-\Delta)^{\frac{-\alpha}{2}} \mathbb{P} ( \f)  \in L^\infty_{0}(\Rt)\cap L^\infty_{4-\alpha}(\Rt)$.  Then, the solution $\u$ has the following asymptotic profile as $|x|\to +\infty$:
\begin{equation}\label{sep 14_3}
    \u(x) = 
    (-\Delta)^{-\frac{\alpha}{2}} \mathbb{P}(\f)(x)
    +
    m_{\alpha} (x) :
    \int_{\Rt} (\u \otimes \u )(y)dy 
  + 
\begin{cases}
     O\left(
\frac{1}{|x|^{9-3\alpha}}
    \right) & \text{if }\alpha \neq 2,
     \\[9pt]
       O\left(
\frac{\log |x|}{|x|^{5-\alpha}}
    \right) & \text{if }\alpha = 2.
\end{cases}
\end{equation}
\end{Proposition}

\medskip 

This  asymptotic profile will  allow us to prove the next theorem.
\begin{Theoreme}\label{Prop-non-existence}
Let $1<\alpha< \frac{5}{2}$. There exists an external force $\f \in \mathcal{S}(\Rt)$, such that the associated solution  $\u$ does not belong to $L^{p,q}(\Rt)$ for all $1\leq p \leq \frac{3}{4-\alpha}$ and $1\leq q <+\infty$. Moreover, for $1\ll |x|$ this solution verifies the estimate from below $ \frac{1}{|x|^{4-\alpha}}   \lesssim 
|\u (x)|$.  

\end{Theoreme} 
\medskip

We continue with the qualitative study of weak $L^p$-solutions to equations (\ref{Frac-NS}). To introduce our next result, first it is worth recalling that in the setting of Sobolev spaces  regularity of weak $\dot{H}^{\frac{\alpha}{2}}$-solutions has been studied in \cite{tang2016new}. Precisely, in the spirit of the celebrated Caffarelli,  Kohn and  Nirenberg theory \cite{cknpseudo}, the authors show a \emph{partial} H\"older regularity result for \emph{suitable} weak  $\dot{H}^{\frac{\alpha}{2}}$-solutions, with $1<\alpha <2$. On the other hand, very recently the authors of \cite{chamorropoggi2023p} show \emph{global} regularity of weak  $\dot{H}^{\frac{\alpha}{2}}$-solutions to equation (\ref{Frac-NS}) in the \emph{homogeneous case} when $\f \equiv 0$. Specifically, the range of values $\frac 5 3<\alpha < 2$ and the fact that $\f \equiv 0$ allow them  to apply a bootstrap argument and to conclude that weak $\dot{H}^{\frac{\alpha}{2}}$-solutions are $C^\infty$-solutions. This argument also  holds in the range $1<\alpha \leq \frac{5}{3}$ as long as the (strong) supplementary hypothesis $\u \in L^{\infty}(\Rt)$ is assumed.  
Similar results on partial regularity of suitable weak solutions for this equation 
in dimension four and five has been studied in \cite{MR4207438,MR3723393}.

\medskip

In our next result, we shall perform  a different approach to prove  \emph{global} regularity of weak $L^p$-solutions to fractional stationary Navier-Stokes equations. We emphasize that this theorem is independent of the previous ones, since here we shall \emph{assume} the existence of weak $L^p$-solutions and we are mainly interested in studying their maximum gain of regularity from an initial regularity on the external force.  Of course, this result holds for solutions constructed in Theorem \ref{Th1}.

\begin{Theoreme}\label{Th3} Let $1<\alpha$ and let $\max\left(\frac{3}{\alpha-1},\frac{6}{\alpha},1 \right)<p<+\infty$. Let $\f \in \dot{W}^{-1,p}(\Rt)$ be an external force and let $\u \in L^p(\Rt)$ be a weak solution of equation (\ref{Frac-NS}) associated to $\f$.  If the external force verifies  $\f \in \dot{W}^{s,p}(\Rt)$, with $0\leq s$, then it holds 
$\u \in \dot{W}^{s+\alpha,p}(\Rt)$ and $P \in \dot{W}^{s+\alpha,p}(\Rt) +  \dot{W}^{s+1,p}(\Rt)$.   
\end{Theoreme} 

Let us briefly explain the general strategy of the proof. The information $\u \in L^p(\Rt)$, with the assumption $\max\left(\frac{3}{\alpha-1},\frac{6}{\alpha},1\right)<p<+\infty$, and the framework of parabolic (time-dependent) fractional Navier-Stokes equations, allows us to prove that $\u \in L^\infty(\Rt)$. Consequently, we are able to get rid of this supplementary hypothesis used in \cite{chamorropoggi2023p}. Thereafter, we use a sharp bootstrap argument  applied to  equations (\ref{Frac-NS-fixed-point}) and (\ref{Pressure}) to obtain the regularity stated above. 

\medskip 

In this regularity result, the assumption $\f \in \dot{W}^{-1,p}(\Rt)$ is essentially technical. On the other hand, the parameter $0\leq s$ measures the initial regularity of the external force $\f$, which yields a gain of regularity of weak $L^p$-solutions to the order $s+\alpha$. This (expected)  \emph{maximum} gain of regularity is given by the effects of the fractional Laplacian operator in equation (\ref{Frac-NS}). See Remark \ref{Rmk-regularity} below for all the details.  Moreover, as pointed out in \cite{chamorropoggi2023p,tang2016new}, the study of regularity in the case $0<\alpha \leq 1$ is (to our knowledge) a hard open problem. 
 
\medskip

When handling the fractional case (with $\alpha \neq 2$) one of the key tools to prove this result is the \emph{fractional Leibniz rule}, so-called the \emph{Kato-Ponce inequality}, which is stated in Lemma \ref{Leibniz-rule} below. To our knowledge, this inequality is unknown for larger spaces than the $L^p$-ones and this fact imposes our hypothesis  $\u \in L^p(\Rt)$.  However, at Appendix \ref{Appendix}, we show that  the classical case (with $\alpha=0$) allows us to prove a more general regularity result for  equation (\ref{Classical-NS}) in the larger  setting of Morrey spaces. 

\medskip

As a direct corollary of Theorem \ref{Th3}, for the particular homogeneous case when $\f \equiv 0$: 
\begin{equation}\label{Frac-NS-homog}
(-\Delta)^{\frac{\alpha}{2}} \u + (\u \cdot \vec{\nabla}) \u + \vec{\nabla}P=0, \quad \text{div}(\u)=0,  \quad \alpha >1,  
\end{equation}
and for the range of values $\max\left(\frac{3}{\alpha-1},\frac{6}{\alpha},1\right)<p<+\infty$, we obtain that weak $L^p$-solutions to this equation are $\mathcal{C}^{\infty}$-functions. This particular result is of interest in connection to another important problem related to equation (\ref{Frac-NS-homog}). 

\medskip

Regularity of weak solutions is one of the key assumptions in the study of their uniqueness. We easily  observe that $\u =0$ and $P=0$ is a trivial solution to (\ref{Frac-NS-homog}) and we look for some functional spaces in which this solution in the unique one. This problem is  so-called the Liouville-type problem for \emph{fractional} stationary Navier-Stokes equations. A \emph{formal} energy estimate and the divergence-free property of velocity $\u$ \emph{heuristically} show that $\| \u \|_{\dot{H}^{\frac{\alpha}{2}}}=0$, which yields to conjecture that the only $\dot{H}^{\frac{\alpha}{2}}$-solution to this equation is the trivial one $\u =0$. However, this problem  is still  out of reach even in the classical case when $\alpha=2$. 

\medskip

We thus look for some \emph{additional hypothesis} on the velocity $\u$ which yield the wished identity $\u=0$. For the \emph{classical}  stationary Navier-Stokes equations (with $\alpha=2$) there is a large amount of literature on the Liouville  problem, see for instance  \cite{chamorrojarrinlem2018p,jarrin2020p,jarrin2023p,seregin2016new} and the references therein. However, Liouville-type problems for the \emph{fractional} case (with $\alpha \neq 2$)  have  been much less  studied by the  research community. In \cite{wang2018p}, the authors show that smooth solutions to (\ref{Frac-NS-homog}) vanish identically when $\u \in \dot{H}^{\frac{\alpha}{2}}(\Rt) \cap L^{\frac{9}{2}}(\Rt)$, with $0<\alpha <2$. Recently, for the range $\frac{5}{3}\leq \alpha <2$, in \cite{yang2022liouville} it is proven that smooth solutions to the equation (\ref{Frac-NS-homog}) are equal to the trivial one, only provided that they belong to the space $L^{\frac{9}{2}}(\Rt)$. As an interesting corollary, due to the Sobolev embedding $\dot{H}^{\frac{5}{6}}(\Rt)\subset L^{\frac{9}{2}}(\Rt)$, for $\alpha=\frac{5}{3}$ the Liouville problem of equation (\ref{Frac-NS-homog}) is solved in the natural energy space $\dot{H}^{\frac{5}{6}}(\Rt)$.  

\medskip 

On the other hand,  for the range of values $\frac{3}{5}<\alpha <2$ the authors of  \cite{chamorropoggi2023p} solve the Liouville-type problem when $\u \in \dot{H}^{\frac{\alpha}{2}}(\Rt) \cap L^{p(\alpha)}(\Rt)$. Here,  the parameter $p(\alpha)$ depends on $\alpha$ and it is close (in some sense) to the critical value $\frac{6}{3-\alpha}$. By Sobolev embeddings  we have $\dot{H}^{\frac{\alpha}{2}}(\Rt) \subset L^{\frac{6}{3-\alpha}}(\Rt)$ and this last value is the natural one to solve this problem. However, the identity $p(\alpha)=\frac{6}{3-\alpha}$ cannot be reached and the information $\u \in L^{p(\alpha)}(\Rt)$ remains an additional hypothesis. 

\medskip

In our next result, we study the Liouville-type problem for a different range of values of $p$. 
\begin{Proposition}\label{Proposition-Liouville} Let $\frac{5}{3}<\alpha \leq 2$, let $\max\left(\frac{3}{\alpha-1},\frac{6}{\alpha}\right)<p\leq \frac{9}{2}$ and let $\u \in \dot{H}^{\frac{\alpha}{2}}(\Rt) \cap L^p(\Rt)$ be a solution to the equation (\ref{Frac-NS-homog}). Then it holds $\u=0$ and $P=0$. 
\end{Proposition}

The proof is mainly based on some \emph{Caccioppoli}-type estimates established in \cite{chamorrojarrinlem2018p} and \cite{chamorropoggi2023p}. These estimates work for smooth enough solutions and in this sense Theorem \ref{Th3} is useful. On the other hand, for the classical case $\alpha=2$ we recover the result proven in \cite[Theorem $1$]{chamorrojarrinlem2018p} (for $3<p\leq \frac{9}{2}$) and this proposition can be seen as its generalization to the fractional case. It is worth mentioning the constraint $\frac{5}{3}<\alpha$ ensures that $\frac{3}{\alpha-1}<\frac{9}{2}$ and $\frac{6}{\alpha}<\frac{9}{2}$. Finally,  by following some of the ideas in \cite{jarrin2023p} we think that this result can be extended for $\frac{9}{2}<p$. 

\medskip

\noindent{\bf Organization of the rest of the paper}. 
The Section \ref{section2} is essentially devoted to present the necessary preliminaries to deal with the proofs of our main results.
In Section \ref{section3}, we prove the existence of weak solutions solutions of the stationary fractional Navier-Stokes in the setting of Lebesgue and Lorentz spaces stated in Theorem \ref{Th1}.
In Section \ref{section4}, we present the proofs of the pointwise estimates and asymptotic profiles stated in Theorem \ref{Th2} and Proposition \ref{Proposition-asymptotic-profile}, respectively.
In Section \ref{section5}, we prove the nonexistence result stated in Theorem \ref{Prop-non-existence}.
While in Sections \ref{section6} and \ref{section7} we prove the Theorem \ref{Th3} and Proposition \ref{Proposition-Liouville}, related to the regularity of weak solutions and a Liouville-type result, respectively.
At the of the paper we present an appendix where we state a new regularity criterion for the classical stationary Navier-Stokes equations in the setting of Morrey spaces.

\medskip 

\paragraph{\bf Acknowledgements.} 
	The second author is supported by the ANID postdoctoral program BCH 2022 grant No. 74220003. Moreover, the authors thank the anonymous referee for  the careful reading of this work and for  accurate corrections   that significantly  improved the article.

\medskip

\paragraph{\bf Statements and Declarations.} Data sharing does not apply to this article as no datasets were generated or analyzed during the current
study. 
In addition, the authors declare that they have no
conflicts of interest and all of them have equally contributed to this paper.

\section{Preliminaries}\label{section2}
In this section,  we summarize some well-known results which will be useful in the sequel.  We start by the Young inequalities in the Lorentz spaces. In particular, the involved constants are \emph{explicitly} written since they will play a substantial role in our study. For a proof we refer to \cite[Section 1.4.3]{chamorro2020pseudo}.

\begin{Proposition}[Young inequalities]\label{Prop-Young} Let $1<p,p_1,p_2<+\infty$ and $1\leq q,q_1, q_2 \leq+\infty$. There exists a generic constant $C>0$ such that  the following estimates hold:  
\begin{enumerate}
\item $\ds{\Vert f \ast g \Vert_{L^{p,q}} \leq C_1 \Vert f \Vert_{L^{p_1, q_1}}\, \Vert g \Vert_{L^{p_2, q_2}}}$, with $1+\frac{1}{p}=\frac{1}{p_1}+\frac{1}{p_2}$, $\frac{1}{q}\leq \frac{1}{q_1}+\frac{1}{q_2}$,  and   $C_1=C\,p\left( \frac{p_1}{p_1-1}\right)\left( \frac{p_2}{p_2-1}\right)$.  
\item $\ds{\Vert f \ast g \Vert_{L^{p,q}} \leq C_2 \Vert f \Vert_{L^1} \, \Vert g \Vert_{L^{p,q}}}$, with $C_2= C\, \frac{p^2}{p-1}$.  
\item $\ds{\Vert f \ast g \Vert_{L^\infty} \leq C_3 \Vert f \Vert_{L^{p,q}}\, \Vert g \Vert_{L^{p',q'}}}$, with  $1=\frac{1}{p}+\frac{1}{p'}$, $1 \leq \frac{1}{q}+\frac{1}{q'}$ and $C_3= C \left(\frac{p}{p-1} \right)\left(\frac{p'}{p' -1} \right)$.
\end{enumerate} 
\end{Proposition} 

Next, the key result in the proof of our theorems related to the existence of solutions is the classical  Picard fixed point scheme, which  we present in the following

\begin{Theoreme}[Picard's fixed point]\label{Th-Picard} Let $(E, \Vert \cdot \Vert_E)$ be a Banach space and let $u_0 \in E$ be an initial data such that $\Vert u_0 \Vert_{E} \leq \delta$. Assume that $B: E \times E \to E$ is a bilinear application. Assume moreover that 
\[ \Vert B(u,v) \Vert_{E} \leq C_B \Vert u \Vert_{E}\Vert u \Vert_{E},\]
for all $u,v \in E$. If there holds 
\[ 0< 4 \delta C_B <1, \]
then the equation $\u=B(\u,\u)+\u_0$ admits a  solution $\u \in E$, which is the unique solution verifying $\Vert \u \Vert_{E} \leq 2\delta$. 
\end{Theoreme}

Thereafter, when studying the regularity of weak solutions we shall use the following lemmas.

\begin{Lemme}[Lebesgue-Besov  Embedding]\label{Lem-U2} Let $1<p < +\infty$. Then, for all $t>0$,  there exist $C>0$ such that the following estimate holds  $$\ds{t^{\frac{3}{\alpha p}} \Vert e^{-t(-\Delta)^{\frac{\alpha}{2}}}  f \Vert_{L^{\infty}} \leq C \Vert f \Vert_{L^{p}}}.$$
\end{Lemme} 	  
\noindent\pv

\noindent Recall first the continuous embedding $L^{p}(\Rt) \subset \dot{B}^{-\frac{3}{p}, \infty}_{\infty}(\Rt)$ (see for instance \cite[Page $171$]{pglemarie2016pseudo}) where the homogeneous Besov space $\dot{B}^{-\frac{3}{p}, \infty}_{\infty}(\Rt)$ can be characterized  as the space of temperate distributions $f\in \mathcal{S}'(\Rt)$ such that $\ds{\| f \|_{\dot{B}^{-\frac{3}{p}, \infty}_{\infty}}= \sup_{t>0}\, t^{\frac{3}{2 p}}\Vert e^{t \Delta}f \Vert_{L^{\infty}} <+\infty}$. Thereafter, by \cite[Page $9$]{pglemarie2018p} we have equivalence $\ds{\| f \|_{\dot{B}^{-\frac{3}{p}, \infty}_{\infty}} \simeq  \sup_{t>0}\, t^{\frac{3}{\alpha p}}\Vert e^{-t(-\Delta)^{\frac{\alpha}{2}}}  f \Vert_{L^{\infty}}}$, from which we  obtain the wished estimate. \finpv 

\medskip  


\medskip

\begin{Lemme}[Fractional Leibniz rule]\label{Leibniz-rule} Let $\alpha>0$, $1<p<+\infty$ and $1<p_0, p_1, q_0, q_1 \leq +\infty$. Then, there exist $C>0$ such that the following estimate holds
\begin{equation*}
\| (-\Delta)^{\frac{\alpha}{2}}(fg) \|_{L^p} \leq C\, \| (-\Delta)^{\frac{\alpha}{2}} f \|_{L^{p_0}}\, \| g \|_{L^{p_1}}+ \| f \|_{L^{q_0}}\, \| (-\Delta)^{\frac{\alpha}{2}} g \|_{L^{q_1}},    
\end{equation*}
where $\frac{1}{p}=\frac{1}{p_0}+\frac{1}{p_1}=\frac{1}{q_0}+\frac{1}{q_1}$. 
\end{Lemme} 
The proof of this estimate can be consulted in \cite{grafakos2014p} or \cite{naibo2019p}.
Finally, given an $\alpha>0$, let $p_\alpha(t,x)$  the convolution kernel of the operator $e^{-t (-\Delta)^{\frac{\alpha}{2}}}$ and $K_\alpha(t,x)=(K_{\alpha, i,j,k}(t,x))_{1\leq i,j,k\leq 3}$  the tensor of convolution kernels associated to the operator $e^{-t (-\Delta)^{\frac{\alpha}{2}}} \P(div(\cdot))$. 

\begin{Lemme}[Lemma $2.2$ of \cite{yu2012p}]\label{Lem-Kernel}  For all $t>0$, there exist a numerical constant $C>0$ depending on  $\alpha$  such that the following estimates hold:
\begin{enumerate}
    \item $\ds{\left\| p_\alpha(t,\cdot) \right\|_{L^1} \leq C}$,
    \item $\ds{\left\| \vec{\nabla} p_\alpha(t,\cdot) \right\|_{L^1} \leq C\,  t^{-\frac{1}{\alpha}}}$, 
    \item  $\ds{\left\|  K_\alpha(t,\cdot) \right\|_{L^1} \leq C\, t^{-\frac{1}{\alpha}}}$.
\end{enumerate}
\end{Lemme}  

\section{Existence of weak solutions in  Lebesgue and Lorentz spaces: proof of Theorem \ref{Th1}}\label{section3}
First recall  that  equation (\ref{Frac-NS}) can be rewritten as the (equivalent) fixed point problem (\ref{Frac-NS-fixed-point})   where,  for the sake of simplicity, we shall denote 
\begin{equation}\label{Bilinear-term}
B(\u,\u)= -(-\Delta)^{-\frac{\alpha}{2}}  \P (\text{div}(\u \otimes \u)).
\end{equation} 

In the next propositions, we shall prove each point stated in Theorem \ref{Th1}.  in the next propositions. 
\begin{Proposition}\label{Prop-Existence-Lorentz} Let $1<\alpha < \frac{5}{2}$. Assume that $(-\Delta)^{-\frac{\alpha}{2}} \f \in L^{\frac{3}{\alpha-1}}(\Rt)$. There exists a universal  quantity $\eta_0(\alpha)>0$, which depends on $\alpha$, such that if $\| (-\Delta)^{-\frac{\alpha}{2}} \P(\f)  \|_{L^{\frac{3}{\alpha-1},\infty}}=\delta <\eta_{0}(\alpha)$ then the equation (\ref{Frac-NS-fixed-point}) has a solution $\u \in L^{\frac{3}{\alpha-1},\infty}(\Rt)$ satisfying and uniquely defined by $\| \u \|_{L^{\frac{3}{\alpha-1},\infty}}\leq 2\delta$. 
\end{Proposition}

\noindent \pv

\noindent Let us start by estimating the bilinear term $B(\u,\u)$. For this we shall need the following technical lemma.
\begin{Lemme}\label{Lem-Bilinear} For $1<\alpha < 4$ we have $\ds{B(\u,\u)=m_\alpha \ast (\u \otimes \u)}$, where  $m_\alpha=(m_{i,j,k})_{1\leq i,j,k \leq 3}$ is a tensor with $m_{i,j,k} \in \mathcal{C}^{\infty}\big(\Rt \setminus \{0\}\big) \cap L^{1}_{loc}(\Rt)$  a homogeneous function of degree $\alpha-4$. Moreover, for all $x\neq 0$ we have $| m_\alpha (x) | \leq c | x |^{\alpha-4}$. 
\end{Lemme}
\noindent \pv 

\noindent Recall that we have $\P(\vec{\varphi})= \vec{\varphi} +  (\vec{\mathcal{R}} \otimes  \vec{\mathcal{R}}) \vec{\varphi}$, where $\vec{\mathcal{R}}=(\mathcal{R}_{i})_{1\leq i \leq 3}$ with $\mathcal{R}_i=\frac{\partial_i}{\sqrt{-\Delta}}$ the i-th Riesz transform.  Then we obtain
\begin{equation*}
 \begin{split}
  B(\u, \u)= &\,  - \ds{\frac{1}{(-\Delta)^{\frac{\alpha}{2}}} \left(  \P (\text{div}(\u \otimes \u))\right)}= \left( -\sum_{j=1}^{3} \sum_{k=1}^{3} (\delta_{i,j}+\mathcal{R}_i\mathcal{R}_j)\frac{\partial_k}{(-\Delta)^{\frac{\alpha}{2}}}(u_j\, u_k)\right)_{1\leq i \leq 3}\\
  =&\, \left( \sum_{j=1}^{3} \sum_{k=1}^{3} m_{i,j,k}\ast(u_j\, u_k)\right)_{1\leq i \leq 3} = m_\alpha \ast (\u \otimes \u).
 \end{split}
\end{equation*}
By a simple computation we write $\ds{\widehat{m}_{i,j,k}(\xi)= - \left( \delta_{i,j}-\frac{\xi_{i} \xi_{j}}{|\xi|^2} \right) \frac{\textbf{i}\xi_k}{|\xi|^\alpha}}$. We  thus have that $\widehat{m}_{i,j,k}\in \mathcal{C}^{\infty}(\Rt \setminus \{0\})$ is a homogeneous function of degree $1-\alpha$ and then $m_{i,j,k}\in \mathcal{C}^{\infty}(\Rt \setminus \{0\})$ is homogeneous function of degree $\alpha-4$.  

\medskip

With this information, for all $x\neq  0$ we can write $| m_\alpha (x)| \leq c | x |^{\alpha-4}$. Moreover,   since $1<\alpha$ we have $-3<\alpha-4$, which yields $m_\alpha \in L^{1}_{loc}(\Rt)$.  Lemma  \ref{Lem-Bilinear} is proven.  \finpv

\medskip 

Observe that  the estimate  $| m_\alpha (x)| \leq c |x|^{\alpha-4}$ for all $x \neq 0$  allows us to conclude that $m_\alpha \in L^{\frac{3}{4-\alpha},\infty}(\Rt)$. Then, we apply  the first point of Proposition \ref{Prop-Young},  where we set the parameters  $p=\frac{3}{\alpha-1}$, $p_1=\frac{3}{4-\alpha}$ and $p_2=\frac{p}{2}=\frac{3}{2(\alpha-1)}$, and we obtain  
\begin{equation*}
\| m_\alpha \ast (\u \otimes \u) \|_{L^{\frac{3}{\alpha-1},\infty}} \leq  C_1 \| m_\alpha \|_{L^{\frac{3}{4-\alpha},\infty}}\, \| \u \otimes \u \|_{L^{\frac{3}{2(\alpha-1)}}}. 
\end{equation*}
\begin{Remarque}\label{Rmk-Constraint-alpha} The condition $1<p_2$ yields the constraint $\alpha < \frac{5}{2}$. We thus set $1<\alpha < \frac{5}{2}$. 
\end{Remarque}

Thereafter, by the H\"older inequalities and by the identity $B(\u, \u)=m_\alpha \ast (\u \otimes \u)$ we have 
\begin{equation}\label{Continuity-bilinear}
\| B(\u, \u) \|_{L^{\frac{3}{\alpha-1},\infty}} \leq C_B(\alpha) \| \u \|^{2}_{L^{\frac{3}{\alpha-1},\infty}}, \qquad \mbox{with} \ \ C_B(\alpha)= C_1\| m_\alpha \|_{L^{\frac{3}{4-\alpha},\infty}}.    
\end{equation}

Thus, for the constant $C_B(\alpha)$ give above, now we set $\delta$ small enough: 
\begin{equation}\label{Condition-delta-1}
  \delta < \frac{1}{4 C_B(\alpha)}=\eta_0(\eta),\end{equation} 
and  by Theorem \ref{Th-Picard} we finish the proof of Proposition  \ref{Prop-Existence-Lorentz}. \finpv 

\medskip

With this first result, now we are able to prove that this  solution  also belongs to the space $L^p(\Rt)$ with $\frac{3}{4-\alpha}<p<+\infty$. In the (more general) case of the Lorentz space $L^{p,q}(\Rt)$ (with $1\leq q \leq +\infty$) the proof follows the same lines, so it is enough to   focus in the case of the Lebesgue spaces. 

\begin{Proposition}\label{Prop-Existence-Lp} With the same hypothesis of Proposition \ref{Prop-Existence-Lorentz}, assume in addition that $(-\Delta)^{-\frac{\alpha}{2}} \f \in L^p(\Rt)$ with $\frac{3}{4-\alpha}<p<+\infty$. There exists a universal quantity $\varepsilon_{0}(\alpha)<\eta_0(\alpha)$, which only depends on $\alpha$, such that if $\delta <\varepsilon_{0}(\alpha)$ then the solution $\u \in L^{\frac{3}{\alpha-1},\infty}(\Rt)$ to the equation (\ref{Frac-NS-fixed-point}) constructed in Proposition \ref{Prop-Existence-Lorentz} verifies $\u \in L^p(\Rt)$ with $\frac{3}{4-\alpha}<p<+\infty$. 
\end{Proposition}

\noindent \pv  

\noindent The solution $\u \in L^{\frac{3}{\alpha-1},\infty}(\Rt)$  to the problem (\ref{Frac-NS-fixed-point}) is obtained as the limit of  the sequence $(\u_n)_{n \in \mathbb{N}}$, where
\begin{equation}\label{Def-Sequence}
\u_{n+1}= B(\u_n, \u_n)+ \u_0, \quad  \mbox{for} \ \  n \geq 0 \ \ \mbox{and} \ \  \u_0= (-\Delta)^{-\frac{\alpha}{2}} \P (\f). 
\end{equation}
We shall use this sequence to prove that $\u \in L^p(\Rt)$. For this we write 
\begin{equation*}
\|  \u_{n+1} \|_{L^p} \leq \| B(\u_n, \u_n) \|_{L^p}+\| \u_0 \|_{L^p}.   
\end{equation*}
For the last expression above, we use our assumption  $(-\Delta)^{-\frac{\alpha}{2}} \f  \in L^p(\Rt)$ to directly obtain $\u_0 \in L^p(\Rt)$. On the other hand, in order to estimate the bilinear term,
we recall that by Lemma \ref{Lem-Bilinear} we have $B(\u_n, \u_n)= m_\alpha \ast (\u_n \otimes \u_n)$ where $m_\alpha \in L^{\frac{3}{4-\alpha}, \infty}(\Rt)$. Therefore, by the first point of Proposition \ref{Prop-Young} (with the parameters $p_1=\frac{3}{4-\alpha}$, $p_2=\frac{3p}{3+p(\alpha-1)}$ and $q=p$,  $q_1=+\infty$, $q_2=p$) we have 
\begin{equation}\label{Def-C1}
\begin{split}
&\| B(\u_n, \u_n) \|_{L^p} \leq C_1(\alpha,p) \| m_\alpha \|_{L^{\frac{3}{4-\alpha},\infty}} \, \| \u_n \otimes \u_n \|_{L^{\frac{3p}{3+p(\alpha-1)},p}}, \\
&\, C_1(\alpha,p)=Cp\left(\frac{3}{7-\alpha}\right)\left(\frac{3p}{(4-\alpha)p-3}\right).    
\end{split}   
\end{equation}
\begin{Remarque} The constant $C_1(\alpha,p)$ defined above blows-up at $p=\frac{3}{4-\alpha}$ and $p=+\infty$. This fact yields the constraint $\frac{3}{4-\alpha}<p<+\infty$. 
\end{Remarque}
Thereafter, from the estimate given in (\ref{Def-C1}) and by the H\"older inequalities we can write 
\begin{equation}\label{eq.bili.lp.lorentz}
\| B(\u_n, \u_n) \|_{L^p} \leq C_1(\alpha,p) \| \u_n  \|_{L^p}\, \| \u_n \|_{L^{\frac{3}{\alpha-1},\infty}}.    
\end{equation}
Moreover,  recall that the sequence defined in (\ref{Def-Sequence}) verifies the uniform estimate $\| \u_n \|_{L^{\frac{3}{\alpha-1},\infty}} \leq 2\delta$, where  $\delta=\left\| (-\Delta)^{-\frac{\alpha}{2}} \P( \f)  \right\|_{L^{\frac{3}{\alpha-1},\infty}}$. Then we get 
\begin{equation*}
\| B(\u_n, \u_n) \|_{L^p} \leq 2\, \delta\, C_1(\alpha,p) \| \u_n  \|_{L^p}.
\end{equation*}
We thus have the following recursive estimate 
\begin{equation}\label{Iterative-un}
\| \u_{n+1} \|_{L^p}  \leq 2\, \delta\, C_1(\alpha,p) \| \u_n  \|_{L^p} + \| \u_0 \|_{L^p}, \ \ \mbox{for all} \ \ n \geq 0.
\end{equation}

At this point, we need to find an additional constraint on the parameter $\delta$  to obtain the control
\begin{equation}\label{Control}
2\, \delta\, C_1(\alpha,p) < \frac{1}{2}.    
\end{equation}
Below we shall need this inequality. Now, in order to get  this control, we need to consider the following cases of the parameter $p \in (\frac{3}{4-\alpha},+\infty)$. First, note  that since $\alpha < \frac{5}{2}$  we have $\frac{3}{4-\alpha} <2<\frac{3}{\alpha-1}$,  and then,  we  split the interval $(\frac{3}{4-\alpha},+\infty)=(\frac{3}{4-\alpha},2) \cup [2, \frac{9}{\alpha-1}] \cup  (\frac{9}{\alpha-1}, +\infty)$.   Moreover, only for technical reasons, first we need to consider the interval $[2, \frac{9}{\alpha-1}]$ since the estimates in the other intervals are based on the ones proven in $[2, \frac{9}{\alpha-1}]$.

\medskip

{\bf The case $p \in [2, \frac{9}{\alpha-1}]$}.  We get back to the expression   of the quantity $C_1(\alpha,p)$ given in (\ref{Def-C1}) and we define  $\ds{0<M(\alpha)=\max_{ p\in [2, \frac{9}{\alpha-1}]} C_1(\alpha,p)<+\infty}$. Then, we set the additional  constraint on the parameter $\delta$ (which already verifies (\ref{Condition-delta-1})) as follows:
\begin{equation}\label{Condition-delta-2}
2 \delta  M(\alpha) < \frac{1}{2},  
\end{equation}
and for all $p\in [2, \frac{9}{\alpha-1}]$ we have the wished control (\ref{Control}).  Therefore, we get back to the inequality (\ref{Iterative-un}) to write 
\begin{equation*}
\| \u_{n+1} \|_{L^p} \leq \frac{1}{2} \| \u_{n} \|_{L^p}+ \| \u_0 \|_{L^p}, \ \ \mbox{for all} \ \ n\geq 0,    
\end{equation*}
hence we obtain the uniform control
\begin{equation}\label{Uniform-control-un}
\| \u_{n+1} \|_{L^p} \leq \left( \sum_{k=0}^{+\infty} \frac{1}{2^k} \right) \, \| \u_0 \|_{L^p}, \ \ \mbox{for all} \ \ n\geq 0.   
\end{equation}
We conclude that the sequence $(\u_n)_{n\in \mathbb{N}}$ is uniformly bounded in the space $L^p(\Rt)$ and then we have $\u \in L^p(\Rt)$. Moreover, recall that this fact holds  as long as the  parameter $\delta$ verifies both conditions (\ref{Condition-delta-1}) and (\ref{Condition-delta-2}), which can be jointly written as
\begin{equation}\label{Condition-delta} 
 \delta < \min\left( \frac{1}{4 C_{B}(\alpha)}, \frac{1}{4 M(\alpha)}\right)=\varepsilon_0(\alpha).    
\end{equation}
The constraint above allowed us to prove that $\u \in L^p(\Rt)$ with $p \in [2, \frac{9}{\alpha-1}]$. Now, we shall use  this same constraint   to prove that $\u \in L^p(\Rt)$ in the  cases $p \in (\frac{3}{4-\alpha},2)$ and $p\in (\frac{9}{\alpha-1}, +\infty)$. 

\medskip

{\bf The case $p \in (\frac{3}{4-\alpha},2)$}.  We get back to equation (\ref{Frac-NS-fixed-point}). Recall that our assumption $(-\Delta)^{-\frac{\alpha}{2}}(\f) \in L^p(\Rt)$ yields that $\u_0 \in L^p(\Rt)$ with $p \in (\frac{3}{4-\alpha},2)$.  On the other hand, recalling that  $(-\Delta)^{-\frac{\alpha}{2}}(\f) \in L^{\frac{3}{\alpha-1},\infty}(\Rt)$  we also  have $\u_0 \in  L^{\frac{3}{\alpha-1},\infty}(\Rt)$  (where $2<\frac{3}{\alpha-1}$) and by a standard interpolation argument we get $\u_0 \in L^2(\Rt)$. Consequently,  by the uniform control given in (\ref{Uniform-control-un}) and by the constraint (\ref{Condition-delta}) we obtain  $\u\in L^2(\Rt)$.  

\medskip

With this information at out disposal, we can prove that $B(\u, \u ) \in L^{\frac{3}{4-\alpha},\infty}(\Rt)$. Indeed, since $\u \in L^2(\Rt)$ we get $\u \otimes \u \in L^1(\Rt)$, and moreover, since $B(\u, \u)=m_\alpha \ast (\u \otimes \u)$ with $m_\alpha \in L^{\frac{3}{4-\alpha},\infty}(\Rt)$, by the second point of Proposition \ref{Prop-Young} we have 
\begin{equation}\label{Bilinear-Lorentz-End-Point}
 B(\u, \u) \in L^{\frac{3}{4-\alpha},\infty}(\Rt).    
\end{equation}

\medskip

Finally, by the estimate (\ref{Continuity-bilinear}) we  also have $B(\u, \u) \in L^{\frac{3}{\alpha-1},\infty}(\Rt)$ (recall that $2<\frac{3}{\alpha-1}$) and  by  well-known interpolation inequalities we get $B(\u, \u)\in L^p(\Rt)$ with  $p\in (\frac{3}{4-\alpha},2)$. Consequently, by the identity (\ref{Frac-NS-fixed-point}) we have $\u \in  L^p(\Rt)$.  

\medskip

{\bf The case $p \in  (\frac{9}{\alpha-1}, +\infty)$}.  We follow similar ideas of the previous case. First remark that we have $\u_0 \in L^p(\Rt)$ with $p\in (\frac{9}{\alpha-1}, +\infty)$ (recall that $\frac{9}{\alpha-1}<p$) and since $\u_0 \in L^{\frac{3}{\alpha-1},\infty}(\Rt)$ by the interpolation inequalities  we obtain $\u_0 \in L^{\frac{4}{\alpha-1}} \cap  L^{\frac{9}{\alpha-1}}(\Rt)$ where $ \frac{4}{\alpha-1}, \frac{9}{\alpha-1} \in [2, \frac{9}{\alpha-1}]$. Thus, always by the uniform control (\ref{Uniform-control-un}) and the constraint (\ref{Condition-delta}) we obtain $\u \in L^{\frac{4}{\alpha-1}} \cap L^{\frac{9}{\alpha-1}}(\Rt)$  and then $\u \in L^{\frac{6}{\alpha-1},2}(\Rt)$. 

\medskip

With this information we can prove that $B(\u, \u) \in L^\infty(\Rt)$. Indeed, since $\u \in L^{\frac{6}{\alpha-1},2}(\Rt)$  we obtain $\u \otimes \u \in L^{\frac{3}{\alpha-1},1}(\Rt)$. Moreover, since $m_\alpha \in L^{\frac{3}{4-\alpha},\infty}(\Rt)$ by the third point of Proposition \ref{Prop-Young} we have 
\begin{equation}\label{Bilinear-Linfty}
B(\u,\u) \in L^{\infty}(\Rt).
\end{equation}
Finally, as we also have $B(\u, \u) \in L^{\frac{3}{\alpha-1},\infty}(\Rt)$ we use again the interpolation inequalities to obtain $B(\u,\u) \in L^p(\Rt)$ with $p \in (\frac{9}{\alpha-1}, +\infty)$, which yields $\u \in L^p(\Rt)$.  

\medskip

Proposition \ref{Prop-Existence-Lp} in now proven. \finpv

\medskip

In order to finish the proof of Theorem \ref{Th1}, with the information obtained in the expressions (\ref{Bilinear-Lorentz-End-Point}) and (\ref{Bilinear-Linfty}) we are able to prove our last proposition. 

\begin{Proposition}\label{4octnprop} With the same hypothesis of Proposition \ref{Prop-Existence-Lorentz}, the following statement holds:
\begin{enumerate}
\item[(A)] If $(-\Delta)^{-\frac{\alpha}{2}}(\f) \in L^{\frac{3}{4-\alpha},\infty}(\Rt)$ then we have $\u \in L^{\frac{3}{4-\alpha},\infty}(\Rt)$.

\medskip

\item[(B)] If $(-\Delta)^{-\frac{\alpha}{2}} \P(\f) \in L^{\infty}(\Rt)$ then we have $\u \in L^{\infty}(\Rt)$.
\end{enumerate}
\end{Proposition}
\noindent \pv 

\noindent The proof is straightforward. For the first point, we just recall that we have  $(-\Delta)^{-\frac{\alpha}{2}}(\f) \in L^{\frac{3}{4-\alpha},\infty}(\Rt) \cap L^{\frac{3}{\alpha-1},\infty}(\Rt)$ and by interpolation we  obtain $\u_0 \in L^2(\Rt)$. We thus  have  the information given in (\ref{Bilinear-Lorentz-End-Point}) which yields $\u \in L^{\frac{3}{4-\alpha},\infty}(\Rt)$. The second point follows the same arguments by assuming that $(-\Delta)^{-\frac{\alpha}{2}}\P(\f) \in L^\infty(\Rt)$. Proposition \ref{4octnprop} is proven.  \finpv

\medskip

Gathering together Propositions \ref{Prop-Existence-Lorentz}, \ref{Prop-Existence-Lp} and \ref{4octnprop}  we conclude the proof of 
Theorem \ref{Th1}. \finpv 
 
\section{Pointwise estimates and asymptotic profiles}\label{section4}
\subsection{Proof of Theorem \ref{Th2}} 
We start this section by proving the following useful lemma.  
\begin{Lemme}\label{Lemma_1pointwise}
Let $1<\alpha < \frac{5}{2}$ and let $\theta_1,\ \theta_2$ two positive constants such that $\alpha-1<\theta_1 + \theta_2<3$. Moreover, let $B(\u, \u)$ the bilinear form defined in (\ref{Bilinear-term}). The following statements holds:
\begin{enumerate}
    \item[(1)] There exists a constant $0<C_1=C_1(\alpha,\theta_1,\theta_2)$ such that 
    \begin{equation*}
    \| 
    B(\u,  \u ) \|_{L^\infty_{1-\alpha+ \theta_1 + \theta_2 }} \leq 
    C_1 
    \| 
    \u \|_{L^\infty_{\theta_1  }}
    \| 
    \u
     \|_{L^\infty_{\theta_2  }}.
     \end{equation*}
     \item[(2)]  There exists a constant $0<C_2=C_2(\alpha)$ such that 
     \begin{equation*}
        \| 
    B(\u,  \u ) \|_{L^\infty_{4-\alpha}} \leq 
    C_2 
\left(
\|\u\|^2_{L^\infty_{\frac{3}{2}}}
+
\|\u\|^2_{L^2}
\right). 
\end{equation*}
\end{enumerate}
\end{Lemme}
\noindent \pv

\noindent The fact that $B(\u,  \u )  = m_\alpha * (\u \otimes  \u)$ and Lemma \ref{Lem-Bilinear} yield the estimate 
\begin{equation*}
    B(\u,  \u )  = 
    \int_{\mathbb{R}^3}
    m_\alpha(x-y) : (\u \otimes  \u)(y)
    dy
    \leq 
    C 
    \| 
    \u \|_{L^\infty_{\theta_1  }}
    \| 
    \u \|_{L^\infty_{\theta_2  }}
    \int_{\mathbb{R}^3}
    \frac{1}{|x-y|^{4-\alpha} |y|^{\theta_1 + \theta_2}  } 
    dy.
\end{equation*}
Is easy to see that the last integral can be bounded by $C_1 |x|^{\alpha - 1 - \theta_1 - \theta_2}$, with 
$0<C_1=C_1(\alpha,\theta_1,\theta_1)$.
Then, we can write
\begin{equation*}
    | B(\u,  \u ) |
     \leq 
     \frac{C_1}{   |x|^{1-\alpha + \theta_1 + \theta_2}  }
     \| 
    \u \|_{ L^\infty_{\theta_1  }  }
    \| 
    \u
     \|_{L^\infty_{\theta_2  }    },
\end{equation*}
and thus the first result follows.

To prove the second part, for $x \neq 0$ we start by splitting  the domain 
$\mathbb{R}^3$ into $\{\frac{|x|}{2} >  |y|\}$ and 
$\{\frac{|x|}{2} \leq  |y|\}$, and then 
\begin{equation*}
\begin{aligned}
    B(\u,  \u )  &= 
    \int_{ \frac{|x|}{2} >  |y|  }
    m_\alpha(x-y) : (\u \otimes  \u)(y)
    dy
    +
    \int_{  \frac{|x|}{2} \leq   |y|  }
    m_\alpha(x-y) : (\u \otimes  \u)(y)
    dy\\
    & = I_1 + I_2.
\end{aligned} 
\end{equation*}
To deal with the integral $I_1$, we note that the integration domain $ \{ | y | \leq \frac{|x|}{2}  \}$ yields $| x-y | \geq | x | - |y | \geq |x| - \frac{| x |}{2} = \frac{| x |}{2} $.  By mixing Lemma \ref{Lem-Bilinear} with the last inequality 
we obtain 
\begin{equation*}
    |I_1| \leq  C
  \int_{   \frac{|x|}{2} >  |y|   }
\frac{1}{|x-y|^{4-\alpha}} |\u (y)|^2
    dy \leq  \frac{C}{|x|^{4-\alpha}}
\int_{  \frac{|x|}{2} >  |y|  }
|\u (y)|^2 dy \leq   \frac{C}{|x|^{4-\alpha}}
    \|\u\|^{2}_{L^2}.
\end{equation*}
For the integral $I_2$, we have
\begin{equation*}
    |I_2| \leq C \|\u \|^2_{ L^\infty_{\frac{3}{2}}} 
    \int_{  \frac{|x|}{2} \leq   |y|  } \frac{1}{|x-y|^{4-\alpha}} 
    \frac{1}{|y|^3} dy.
\end{equation*}
The fact that $\frac{|x|}{2} \leq   |y|$ yields  
\begin{equation*}
    \int_{  \frac{|x|}{2} \leq   |y|  } \frac{1}{|x-y|^{4-\alpha}} 
    \frac{1}{|y|^3} dy
    \leq 
\frac{1}{|x|}
    \int_{  \frac{|x|}{2} \leq   |y|  } \frac{1}{|x-y|^{4-\alpha}} 
    \frac{1}{|y|^2} dy.
\end{equation*}
Is straightforward to see that the last integral can be bounded by $C_2|x|^{\alpha-3}$, with $0<C_2=C_2(\alpha)$, and thus we can conclude 
\begin{equation*}
    |I_2| \leq \frac{C_2}{|x|^{4-\alpha}}\|\u \|^2_{L^\infty_{\frac{3}{2}}}.
\end{equation*}
The lemma is then proved.  \finpv 

\medskip

Now, in the next propositions,  we  prove each point stated in Theorem \ref{Th2}. 

\begin{Proposition}\label{Prop-Dec-1} Let $1<\alpha<5/2$. Assume that $(-\Delta)^{-\frac{\alpha}{2}} \P (\f) \in L^{\infty}_{\alpha-1}(\Rt)$. There exists a universal quantity $0<\eta_1(\alpha)<\varepsilon_0(\alpha)$ such that if $\ds{\|(-\Delta)^{-\frac{\alpha}{2}} \P (\f) \|_{L^{\infty}_{\alpha-1}}=\delta < \eta_1(\alpha)}$  then the  equation (\ref{Frac-NS}) has a solution $\u \in L^{\infty}_{\alpha-1}(\Rt)$ satisfying and uniquely defined by the condition
$\| \u \|_{L^{\infty}_{\alpha}} \leq 2 \delta$. 
\end{Proposition}
\noindent \pv

\noindent From Theorem \ref{Th1}  we have a solution $\u$ of equation  \eqref{Frac-NS} in the space $L^{\frac{3}{\alpha-1},\infty}(\Rt)$. In view of our hypotheses for $(-\Delta)^{\frac{\alpha}{2}} \mathbb{P}( \f) \in L^{\infty}_{\alpha-1}(\Rt)$ to conclude that $\u$ belongs to $L^\infty_{\alpha-1}(\Rt)$, we shall to prove 
\begin{equation*}
 \| 
    B(\u,  \v ) \|_{L^\infty_{\alpha-1 }} \leq 
    C 
    \| 
    \u \|_{L^\infty_{\alpha-1   }}
    \| 
    \v
     \|_{L^\infty_{\alpha-1 }} ,
\end{equation*}
with $C$ independent of $\u$ and $\v$. This fact follows directly  by taking $\theta_1=\theta_2=\alpha-1$
in Lemma \ref{Lemma_1pointwise}. Now, considering 
\begin{equation*}
    \|(-\Delta)^{-\frac{\alpha}{2}} \mathbb{P} f\|_{L^\infty_{\alpha-1}}< \eta_1(\alpha)<\varepsilon_0(\alpha),
\end{equation*}
where $\varepsilon_0(\alpha)$
is the constant of Theorem \ref{Th1}, by applying Picard's fixed-point Theorem \ref{Th-Picard} we obtain the existence and uniqueness of the solution $\u$ in $L^\infty_{\alpha-1}(\Rt)$.  \finpv 

\begin{Remarque}
The solution constructed previously also belongs to $L^{\frac{3}{\alpha-1},\infty}(\Rt)$ since $L^\infty_{\alpha-1}(\Rt) \subset L^{\frac{3}{\alpha-1},\infty}(\Rt)$.
\end{Remarque} 

\begin{Proposition}\label{oct4prop1} With the same hypotheses of Proposition \ref{Prop-Dec-1}, assume that $(-\Delta)^{-\frac{\alpha}{2}}\P(\f) \in L^{\infty}_{\theta}(\Rt)$, with  $0<\theta \leq 4-\alpha$. There exists a universal quantity $\varepsilon_1(\alpha)<\eta_1(\alpha)$ such that if $\delta < \varepsilon_1(\alpha)$ then the solution $\u \in L^{\infty}_{\alpha-1}(\Rt)$ constructed by Proposition \ref{Prop-Dec-1} verifies $\u \in L^{\infty}_{\theta}(\Rt)$. 
    \end{Proposition}
\noindent \pv

\noindent By assumption  we have $(-\Delta)^{-\frac{\alpha}{2}} \mathbb{P}(\f) \in L^\infty_{\theta}(\Rt)$, then  considering  the sequence such that $\u_{n+1} =  \u_0 + B(\u_{n},\u_{n})$  (with $\u_0=(-\Delta)^{-\frac{\alpha}{2}} \mathbb{P}(\f)$) and Lemma \ref{Lemma_1pointwise}
with $\theta_1 = \theta, \ \theta_2 = \alpha-1$ 
, we obtain the estimate
\begin{equation}\label{eq. iterative space E}
    \| 
    B(\u_n,  \u_n ) \|_{L^\infty_{\theta }} \leq 
    C_1(\alpha,\theta)
    \| 
    \u_n \|_{L^\infty_{\alpha -1   }}
    \| 
    \u_n
     \|_{L^\infty_{\theta }}
    \quad \text{for} \quad  0<\theta < 4 - \alpha.
\end{equation}
\begin{Remarque}
The positive constant $C_1(\alpha,\theta)$ depends continuously of $\theta\in (0,4-\alpha)$ and it blows-up at $\theta=0$ and $\theta =4-\alpha$.  
\end{Remarque}
On the other hand,  recall that by the Picard iterative schema  applied to the approximation $\u_n$ we have the uniform control
$$
\|\u_n\|_{  L^\infty_{\alpha-1}  } \leq 2 \| \u_0\|_{ L^\infty_{\alpha-1}  }=2\delta. 
$$ 
Then, by  \eqref{eq. iterative space E} we get
\begin{equation}\label{ineq. iterative}
    \| \u_{n+1} \|_{L^\infty_{\theta}}  \leq 2\, \delta\, C_1(\alpha,\theta) \| \u_n  \|_{L^\infty_{\theta}} + \| \u_0 \|_{L^\infty_{\theta}},
    \quad
    \text{for}
    \quad 
    0<\theta < 4 - \alpha.
\end{equation}
As in the proof of Proposition \ref{Prop-Existence-Lp}, we need to find an additional constraint on $\delta$ to get
\begin{equation}\label{control-delta-2}
2\, \delta\, C_1(\alpha,\theta) < \frac{1}{2},
\end{equation}
which we will use later. To obtain this inequality, we split $(0,4-\alpha]=(0, \frac{\alpha-1}{2})\cup [\frac{\alpha-1}{2}, \frac{3}{2}]\cup (\frac{3}{2}, 4-\alpha]$ (remark that $1<\alpha<5/2$ yields $\frac{3}{2}<4-\alpha$). Moreover, only for technical reasons, first we shall consider the case $\theta \in [\frac{\alpha-1}{2}, \frac{3}{2}]$ and then we will study the cases $\theta \in (0, \frac{\alpha-1}{2})$ and $\theta \in (\frac{3}{2}, 4-\alpha]$. 

\medskip 

  {\bf The case $\theta \in [\frac{\alpha-1}{2}, \frac{3}{2}]$}. We define the quantity $\ds{0<N(\alpha)= \max_{\theta \in [\frac{\alpha-1}{2}, \frac{3}{2}]} C_1(\alpha,\theta)<+\infty}$. Then, we set the additional constraint on $\delta$
\begin{equation*}
2\delta N(\alpha)<\frac{1}{2},
\end{equation*}
which yields (\ref{control-delta-2}),
and by following the same arguments in (\ref{Uniform-control-un}) we obtain that $\u \in L^{\infty}_{\theta}(\Rt)$. 

\medskip

  {\bf The case $\theta \in (\frac{3}{2}, 4-\alpha]$}. Recall that $(-\Delta)^{-\frac{\alpha}{2}}\P(\f) \in L^{\infty}_{\alpha-1}(\Rt) \cap L^{\infty}_{\theta}(\Rt)$ and then $(-\Delta)^{-\frac{\alpha}{2}}\P(\f) \in   L^{\infty}_{\frac{3}{2}}(\Rt)$.   
 Thus, by the case above we get $\u \in L^{\infty}_{\frac{3}{2}}(\Rt)$. 

\medskip

  Before going further let consider the following useful result.
\begin{Lemme}\label{aug-lema-1}
    Let $\theta_1, \theta_2>0$ and $p>1$ such that 
    $\frac{3}{\theta_1} 
    <  p  <
    \frac{3}{\theta_2}$. Then 
    $
        L^\infty_{\theta_1}(\Rt)\cap L^\infty_{\theta_2}(\Rt) \subset L^{p}(\Rt). 
    $
\end{Lemme}
\noindent \pv

\noindent The inclusions $L^\infty_{\theta_1}(\Rt) \subset L^{\frac{3}{\theta_1},\infty}(\Rt)$ and $ L^\infty_{\theta_2}(\Rt) \subset L^{\frac{3}{\theta_2},\infty}(\Rt)$
yield
\begin{equation}\label{aug-2-2023}
    L^\infty_{\theta_1}\cap L^\infty_{\theta_2} (\Rt)
    \subset
    L^{\frac{3}{\theta_1},\infty}
    \cap
    L^{\frac{3}{\theta_2},\infty}(\Rt).
\end{equation}
By hypothesis and  interpolation of Lorentz spaces, we can write
\begin{equation}\label{aug-1-2023}
    L^{\frac{3}{\theta_1},\infty}
    \cap
    L^{\frac{3}{\theta_2},\infty}(\Rt)
    \subset 
    L^{p,p}(\Rt).
\end{equation}
Then, by mixing \eqref{aug-2-2023} and \eqref{aug-1-2023} we conclude the proof.
\finpv

\medskip

Thus, considering Lemma \ref{aug-lema-1} and the fact that $(-\Delta)^{-\frac{\alpha}{2}}\P(\f) \in L^{\infty}_{\alpha-1}(\Rt) \cap L^{\infty}_{\theta}(\Rt)$, we get   $(-\Delta)^{-\frac{\alpha}{2}}\P(\f) \in L^2(\Rt)$. Thus,  by Proposition \ref{Prop-Existence-Lp} we have $\u \in L^2(\Rt)$, and by Lemma \ref{Lemma_1pointwise} we obtain $B(\u,\u)\in L^{\infty}_{4-\alpha}(\Rt)$. Recalling that we also have $ B(\u,\u)\in L^{\infty}_{\alpha-1}(\Rt)$, and moreover, since  $\alpha-1<\frac{3}{2}$, we get $B(\u,\u)\in  L^{\infty}_{\theta}(\Rt)$, which yields $\u \in L^{\infty}_{\theta}(\Rt)$. 

\medskip 

  {\bf The case $\theta \in (0, \frac{\alpha-1}{2})$}. Similar arguments as above yield $\u \in L^\infty_{\alpha-1} \cap L^\infty_{\frac{\alpha-1}{2}}$, and then 
$\u \in L^\infty_{\frac{\theta+\alpha-1}{2}}(\Rt)$. From Lemma \ref{Lemma_1pointwise} we obtain $B(\u,\u)\in L^\infty_\theta(\Rt)$, and then $\u\in L^\infty_\theta(\Rt)$. With this we conclude the proof of Proposition \ref{oct4prop1}. \finpv

 \medskip

\noindent Gathering together Propositions \ref{Prop-Dec-1} and \ref{oct4prop1} we finish with  the proof of Theorem \ref{Th2}. \finpv

 \medskip

\noindent In this point we stress the fact that Corollary \ref{Corolary-pointwise-estimate} follows directly as a consequence of Proposition \ref{oct4prop1}. 

 \medskip
\subsection{Proof of Corollary \ref{Corolary-pointwise-estimate}} 
We begin by stressing that there exist a numerical constant $C$ such that
\begin{equation*}
    |(-\Delta)^{\frac{-\alpha}{2}} \P(\f) | \leq {C}
    ,\qquad  
        |(-\Delta)^{-\frac{\alpha}{2}} \P(\f) | \leq \frac{C}{ |x|^{4-\alpha} },
\end{equation*}
 and
\begin{equation*}
    |B(\u,\u)| \leq {C}
    ,\qquad 
        |B(\u,\u) | \leq \frac{C}{ |x|^{4-\alpha} }.
\end{equation*}
Gathering the expressions above with the fact that $\u= (-\Delta)^{-\frac{\alpha}{2}} \P(\f) +
B(\u,\u)
$, we conclude the pointwise estimate 
$$ ( 1  + |x|^{4-\alpha } ) |\u(x)| \leq C ,$$ 
and then we obtain \eqref{sep_1_gaston}. With this we conclude the proof of Corollary \ref{Corolary-pointwise-estimate}. \finpv

\medskip

 \subsection{Proof of Proposition \ref{Proposition-asymptotic-profile}}

With the information obtained in Propositions \ref{Prop-Dec-1} and \ref{oct4prop1}, we are able to derive sharp asymptotic profiles of solutions. 

\medskip

\noindent  \pv

\noindent In the next we will prove that the bilinear term $B(\u,\u)$ has the following asymptotic profile as $|x|\to +\infty$: 
\begin{equation}\label{sep 14_2}
   B(\u,\u) =
    m_{\alpha} (x) :
    \int_{\Rt} (\u \otimes \u )(y)dy 
  + 
\begin{cases}
     O\left(
\frac{1}{|x|^{9-3\alpha}}
    \right) & \text{if }\alpha \neq 2,
     \\[9pt]
       O\left(
\frac{\log |x|}{|x|^{5-\alpha}}
    \right) & \text{if }\alpha = 2.
\end{cases}
\end{equation}
To this end, we consider the following descomposition of the bilinear term $B(\u,\u)$: 
\begin{equation*}
    \begin{split}
  (-\Delta)^{-\frac{\alpha}{2}} \P div (u \otimes u)(x) =&\, 
  \int_{\Rt} m_\alpha(x-y) : (\u \otimes \u )(y) dy \\
  =&\, m_\alpha(x) : \int_{\Rt} (\u \otimes \u )(y) dy \\
  & \, - m_\alpha(x): \int_{|y|\geq |x|/2} (\u \otimes \u )(y) dy \\
  &\, + \int_{|y|\leq |x|/2} (m_\alpha(x-y)-m_\alpha(x)):(\u \otimes \u )(y) dy\\
  &\, + \int_{|x-y|\leq |x|/2} m_\alpha(x-y): (\u \otimes \u )(y) dy, \\
  &\, + \int_{|y|\geq |x|/2, |x-y|\geq |x|/2} m_\alpha(x-y):(\u \otimes \u )(y) dy\\
  = &\,  I_1 + I_2 + I_3 + I_4 + I_5, 
    \end{split}
\end{equation*}
where we must estimate the terms from $I_2$ to $I_5$.
\begin{itemize}
    \item \underline{Term $I_2$.} In this case we have
    \begin{equation}\label{I2f}
        I_2 \leq C |x|^{-(9-3\alpha)},   \quad \text{as } |x|\to +\infty.  
    \end{equation}
    In fact, since $\u\in L^\infty_{4-\alpha}(\Rt)$, we have  $|\u \otimes \u(y)| \leq C|y|^{-2(4-\alpha)}$. Then, considering that $|m_\alpha (x)|\leq C|x|^{-(4-\alpha)}$,  the change of variables $\rho=|y|$ yields 
  \begin{equation*}
\begin{split}
    I_2 \leq &\, |x|^{-(4-\alpha)} \int_{|y|\geq |x|/2} |y|^{-2(4-\alpha)} dy \\
    \leq &\, |x|^{-(4-\alpha)} \int_{|y|\geq |x|/2} |y|^{ -2(4-\alpha)} dy\\
    \leq &\, |x|^{-(4-\alpha)} \int_{|y|\geq |x|/2} |y|^{ -2(4-\alpha)} dy\\
    \leq &\, |x|^{-(4-\alpha)} \int_{|x|/2}^{+\infty} \frac{\rho^2 d \rho}{\rho^{{2(4-\alpha)}}}\\
    \leq &\, C |x|^{-(4-\alpha)} |x|^{-(5-2\alpha)}=C|x|^{-(9-3\alpha)}.
\end{split}
\end{equation*}
    \item \underline{Term $I_3$.} In this case we have 
    \begin{equation}\label{I_3f}
    I_3\leq
\begin{cases}
     \dfrac{ C }{ |x|^{9-3\alpha} } & \text{if }\alpha \neq 2,
     \\[9pt]
      C \dfrac{\log(|x|)}{ |x|^{ 5-\alpha }   } & \text{if }\alpha = 2,
\end{cases}
    \quad \text{as } |x|\to +\infty.
\end{equation}
In fact, since $|\nabla m_\alpha(x)|\leq C |x|^{-(5-\alpha)}$, by the main value theorem with  $z=\theta (x-y)+(1-\theta)x$ (where $0<\theta<1$)  we can write 
\begin{equation*}
     |m_\alpha(x-y)-m_\alpha(x)| \leq C |\nabla  m_\alpha(z)|| (x-y)-x| \leq C |z|^{- (5-\alpha) } |y|.
\end{equation*}
By mixing $z=\theta (x-y)+(1-\theta)x=x-\theta y$ with  $0<\theta<1$ and $|y|\leq |x|/2$, we obtain  
\begin{equation*}
  |z|=| x-\theta y | \geq |x| - \theta |y| \geq |x|-|y|\geq |x|-|x|/2=|x|/2, 
\end{equation*}
 and then 
\begin{equation*}
    |m_\alpha(x-y)-m_\alpha(x)| \leq C |z|^{-(5-\alpha)} |y| \leq C |x|^{-(5-\alpha)}|y|. 
\end{equation*}
Thus, we get 
\begin{equation}\label{11sep3}
    I_3 \leq C|x|^{-(5-\alpha)} \int_{|y|\leq |x|/2} |y| |\u (y)|^2 dy.
\end{equation}
Since  $1 \ll |x|$, we can write 
\begin{equation*}
\int_{|y|\leq |x|/2} |y| |\u (y)|^2 dy= \int_{|y|\leq 1} |y| |\u (y)|^2 dy+\int_{1<|y|\leq |x|/2} |y| |\u (y)|^2 dy=I_{3,1}+I_{3,2}   . 
\end{equation*}
To control the term $I_{3,1}$ we stress the fact that 
$u \in L^2(\Rt)$, then 
\begin{equation}\label{11sep2}
  I_{3,1}= \int_{|y|\leq 1} |y| |\u (y)|^2 \leq \int_{|y|\leq 1}  |\u (y)|^2 \leq \| \u \|^{2}_{L^2}.  
\end{equation}
Now, to deal with $I_{3,2}$, the fact that $|u(y)|\leq C |y|^{-(4-\alpha)}$  and the change of variables $\rho=|y|$ yield 
\begin{equation}\label{11sep1}
     I_{3,2} \leq \int_{1<|y|\leq |x|/2} |y| |\u (y)|^2 dy\leq \int_{1}^{|x|/2} \frac{\rho^2 d \rho}{\rho^{ 7-2\alpha   }}
     = 
\begin{cases}
     \dfrac{ C }{ |x|^{4-2\alpha} } & \text{if }\alpha \neq 2
     \\[7pt]
      C \log(|x|) & \text{if }\alpha = 2
\end{cases}
\end{equation}

Gathering together the estimations \eqref{11sep1} and \eqref{11sep2} in \eqref{11sep3}, we obtain \eqref{I_3f}.

\item \underline{Term $I_4$.} In this case the following pointwise estimate follows
\begin{equation}\label{I4}
    I_4 \leq C |x|^{-(9-3\alpha)},  \quad \text{as } |x|\to +\infty.
\end{equation}

In fact, since $|x-y|\leq |x|/2$ we can write 
\[ |y|= |x-(x-y)| \geq |x| - |x-y| \geq |x| - |x|/2 = |x|/2.\]
Then,  considering $|m(x-y)|\leq C|x-y|^{-(4-\alpha)}$ and $|u(y)|^2 \leq C|y|^{-2(4-\alpha)}$, we obtain 
\begin{equation*}
    \begin{split}
        I_4 & \leq \, \int_{|x-y|\leq |x|/2} |x-y|^{-(4-\alpha)}|y|^{-2(4-\alpha)} dy \\
        & \leq\,  |x|^{-(4-\alpha)} \int_{|x-y|\leq |x|/2} |x-y|^{-(4-\alpha)}|y|^{-(4-\alpha)} dy  
        \\
        & \leq\,  |x|^{-(4-\alpha)} \int_{\Rt} |x-y|^{-(4-\alpha)}|y|^{-(4-\alpha)} dy 
        \\
        & \leq \,  C|x|^{-(4-\alpha)}|x|^{-(5-2\alpha)}
        \\
        & = C|x|^{-(9-3\alpha)}. 
    \end{split}
\end{equation*}

\item \underline{Term $I_5$.} In this case the following pointwise estimate follows
\begin{equation}\label{I5f}
    I_5 \leq C |x|^{-(9-3\alpha)},  \quad \text{as } |x|\to +\infty.
\end{equation}
In fact, since $|y|\geq |x|/2$ and $|x-y| \geq |x|/2$ we can write
\begin{equation*}
    I_5 \leq \int_{|y|\geq |x|/2, |x-y| \geq |x|/2} |m_\alpha(x-y)| \cdot  |y|^{-2(4-\alpha)} dy \leq C|x|^{-(4-\alpha)} \int_{|y|\geq |x|/2}|y|^{-2(4-\alpha)}dy.
\end{equation*}
Considering  similar arguments as in the case of  term $I_2$ we conclude  (\ref{I5f}). 
\end{itemize}
Gathering together all the estimations obtained above we deduce the asymptotic profiles \eqref{sep 14_3} and \eqref{sep 14_2}.

With this we conclude the proof of Proposition \ref{Proposition-asymptotic-profile}. \finpv




\section{Nonexistence result: proof of Theorem \ref{Prop-non-existence}}\label{section5} 
Let us explain the general strategy of the proof. The term $ m_{\alpha} (x) :
    \int_{\Rt} (\u \otimes \u )(y)dy$ in the asymptotic profile (\ref{sep 14_3}), the fact that the tensor $m_\alpha$ is a homogeneous functions of order $4-\alpha$ (see Lemma \ref{Lem-Bilinear} below) and a well-prepared external force $\f$,  will allow us to obtain the  estimate from below
    \[ \frac{1}{|x|^{4-\alpha}} \lesssim | \u(x)|, \quad 1\ll |x|, \]
  from which Theorem \ref{Prop-non-existence} directly follows. 

  \medskip
  
  However, one risks that the term $ m_{\alpha} (x) :
    \int_{\Rt} (\u \otimes \u )(y)dy$ vanishes identically, and our starting point is to study when this fact holds. 
\begin{Proposition}\label{11octprop1}
Under the same hypotheses of Proposition \ref{Prop-Dec-1}, assume that  $(-\Delta)^{\frac{-\alpha}{2}} \mathbb{P} ( \f)  \in L^\infty_{0}(\Rt)\cap L^\infty_{4-\alpha}(\Rt)$. Then,  the term  $m_{\alpha} (x) :
    \int_{\Rt} (\u \otimes \u )(y)dy$ vanishes identically on $\Rt$, if and only if, there exist a constant $c\in \R$ such that, for $j,k =1,2,3$ we have
    \begin{equation}\label{Condition-non-existence}
    \int_{\Rt} u_j u_k = c \, \delta_{j,k} =:  \begin{cases}
     c & \text{if }j=k,\\
      0 & \text{otherwise}.
\end{cases}
\end{equation}
\end{Proposition}
\noindent \pv

\noindent 
Note that, by
considering the term 
$ m_{\alpha} (x) :
    \int_{\Rt} (\u \otimes \u )(y)dy$ (already defined in (\ref{Def-two-dots})),
    and taking the Fourier transform on the Kernel  $ m_{\alpha} (x) $, 
  we can write for $i=1,2,3$
\begin{equation}\label{sept 14_1}
    \begin{aligned}
        \sum_{j=1}^{3} \sum_{k=1}^{3} 
\widehat{m}_{i,j,k}(\xi)
\int_{\Rt} u_ju_k
   & = \sum_{j=1}^{3} \sum_{k=1}^{3} 
    - \left( 
    \delta_{i,j}   -   \frac{ \xi_{i} \xi_{j} }{ |\xi|^2 } 
    \right)
    \frac{ \textbf{i}\xi_k }{ |\xi|^\alpha }
    \int_{\Rt} u_j u_k   
    \\
    &
    =
    R_i(\xi) 
    \frac{\textbf{i}}{|\xi|^{\alpha+2}} 
    \int_{\Rt} u_j u_k , 
    \end{aligned} 
\end{equation}
where $\displaystyle R_i(\xi) = -\left(\sum_{j=1}^{3} \sum_{k=1}^{3}  \delta_{i,j} \xi_k  |\xi|^2 
    - \xi_{i} \xi_{j} \xi_k \right)$.
    Thus, the vanishing condition
$\int_{\Rt} u_j u_k = c \delta_{j,k}$ holds,   if and only if $ m_{\alpha} (x) :
    \int_{\Rt} (\u \otimes \u )(y)dy \equiv 0$. Note that, this last expression follows 
 by considering 
the expression above equal to 0 and the following result of L. Brandolese and F. Vigneron in \cite{brandolese2007new}.
\begin{Lemme}
    For any numerical  matrix $A=(a_{j,k})_{1\leq j,k\leq 3}$, let define the family of homogeneous polynomials  
    \begin{equation*}
Q_i(\xi)=  \sum_{j=1}^{3} \sum_{k=1}^{3}  \left(|\xi|^2  
( \delta_{j,k} \xi_k  +  \delta_{i,k}  \xi_j   +  \delta_{i,j}   \xi_k   )
- 5 \xi_i \xi_j \xi_k\right) a_{j, k}, \quad i=1,2,3. 
\end{equation*}
Then, the following assertions are equivalent:
    \begin{enumerate}
        \item The matrix $A$ is proportional to the identity matrix.
        \item $Q_i \equiv 0$ for all indices $i=1,2,3.$
        \item There exist an index $1\leq i \leq 3$ such that $Q_i \equiv 0$.
        \item There exist an index $1\leq i\leq 3$ such that $\partial_i Q_i \equiv 0$.
    \end{enumerate}    
\end{Lemme}
With this we conclude the proof of Proposition \ref{11octprop1}. \finpv  

\medskip

With this proposition at our disposal,   we shall construct a well-prepared external force $\f \in \mathcal{S}(\Rt)$, such that its associated solution to equation (\ref{Frac-NS}) does not verify the condition (\ref{Condition-non-existence}).  
\begin{Proposition}\label{prop 10oct1} 
Let $\vec{\mathfrak{f}}_0$ be a divergence-free vector field satisfying the following assumptions:
\begin{itemize}
    \item $\widehat{\vec{\mathfrak{f}}_0} \in C_0^\infty(\Rt)$,
    \item $0\notin \text{supp}(\widehat {\vec{\mathfrak{f}}_0})$,
    \item the matrix $\displaystyle \left( 
    \int
    \frac{( \widehat {\mathfrak{f}_0})_i(\overline{ \widehat {\mathfrak{f}_0}})_j}{|\xi|^{2\alpha}} d\xi 
    \right)_{i,j} 
    $ is not a scalar multiple of the identity matrix $\mathbb{I}_3$.
\end{itemize}
 Then, there exist $\eta_0 >0$ such that the solution $\u  \in L^{\frac{3}{\alpha-1},\infty}(\Rt) \cap L^2(\Rt)$ of the stationary Navier-Stokes equation (\ref{Frac-NS}) with $\vec{f} = \eta {\vec{\mathfrak{f}}_0}$ and 
$0< \eta \leq \eta_0$ does not verify (\ref{Condition-non-existence}). 
\end{Proposition} 
\noindent \pv 

\noindent  Let consider the positive constant $\epsilon_1=\epsilon_1(\alpha)$ arising in Theorem \ref{Th1} and $\eta_0>0$ small enough satisfying 
\begin{equation*}
    \eta_0 \left\|  (-\Delta)^{-\frac{\alpha}{2}} \mathbb{P} ( {\vec{\mathfrak{f}}_0})  
    \right\|_{L^{ \frac{3}{\alpha-1}, \infty}}
    \leq \epsilon_1. 
\end{equation*}
Thus, under the framework of Theorem \ref{Th1}, we know that
 for $0< \eta \leq \eta_0$ and $\vec{f} =\eta {\vec{\mathfrak{f}}_0}$, there exist a unique solution $\u \in L^{\frac{3}{\alpha-1},\infty} \cap L^2 $ such that $$\| \u \|_{L^{\frac{3}{\alpha-1},\infty}} \leq 2\eta  \| (-\Delta)^{-\frac{\alpha}{2}} \P ( {\vec{\mathfrak{f}}_0} ) \|_{L^{\frac{3}{\alpha-1},\infty}}.$$ 
We set $\u_0 = \eta  (-\Delta)^{-\frac{\alpha}{2}} \P ( {\vec{\mathfrak{f}}_0} ) $.  
By considering $\u$ in the iteration scheme and  \eqref{eq.bili.lp.lorentz} with $p=2$, we can write 
\begin{equation}\label{sept20eq1}
\begin{aligned}
        \| \u - \u_0    \|_{L^2}  
        & =
    \| B(\u,\u)    \|_{L^2}
    \\
    & \leq  
    C \| \u     \|_{L^2} 
    \|  \u   \|_{L^{\frac{3}{\alpha-1},\infty}}
    \\
    & 
    \leq 
 \| \u     \|_{L^2}\,   2  C \eta 
    \|  (-\Delta)^{-\frac{\alpha}{2}} \P ( {\vec{\mathfrak{f}}_0} )     \|_{L^{\frac{3}{\alpha-1},\infty}}.
\end{aligned}
\end{equation}
Thus,
\begin{equation*}
    \begin{aligned}
        \| \u     \|_{L^2}  
        & \leq 
    \| \u_0   \|_{L^2}
    +
     \|  \u -  \u_0   \|_{L^2}
    \\
    &
    \leq 
\eta 
\|  (-\Delta)^{-\frac{\alpha}{2}} \P ( {\vec{\mathfrak{f}}_0} )     \|_{L^2}
    +
   2 C \eta \| \u     \|_{L^2} 
    \|  (-\Delta)^{-\frac{\alpha}{2}} \P ( {\vec{\mathfrak{f}}_0} )     \|_{L^{\frac{3}{\alpha-1},\infty}}.
\end{aligned}
\end{equation*}
By choosing $\eta_0$ such that 
\begin{equation*}
    \eta_0 \leq 
    \frac{1}{   
    4C  \|  (-\Delta)^{-\frac{\alpha}{2}} \P ( {\vec{\mathfrak{f}}_0} )     \|_{L^{\frac{3}{\alpha-1},\infty}}
    }, 
\end{equation*}
we obtain 
\begin{equation}\label{sept20eq3}
       \| \u     \|_{L^2}  
       \leq 
       2 \eta 
          \|    (-\Delta)^{-\frac{\alpha}{2}} \P ( {\vec{\mathfrak{f}}_0} )   \|_{L^2}.  
\end{equation}
Gathering together the previous estimation with 
\eqref{sept20eq1} we conclude 
\begin{equation}\label{sept20eq2}
      \| \u - \u_0    \|_{L^2} 
      \leq 
      4C \eta^2 
      \|  (-\Delta)^{-\frac{\alpha}{2}} \P ( {\vec{\mathfrak{f}}_0} )     \|_{L^{\frac{3}{\alpha-1},\infty}}
       \|    (-\Delta)^{-\frac{\alpha}{2}} \P ( {\vec{\mathfrak{f}}_0} )   \|_{L^2}.  
\end{equation}

\medskip

To continue, note that our third  assumption on ${\vec{\mathfrak{f}}_0}$  yields
\begin{equation*}
    \int    (-\Delta)^{-\frac{\alpha}{2}} \P ( {\vec{\mathfrak{f}}_0} ) \otimes   (-\Delta)^{-\frac{\alpha}{2}} \P ( {\vec{\mathfrak{f}}_0} )
    \neq 
    \alpha \mathbb{I}_3,\quad
    \text{for }
    \alpha\in \R.
\end{equation*}
This fact means that there exist different indices $i$ and $j$ such that either: 
\begin{itemize}
    \item[$i)$]  $\displaystyle  \int    
\left|
\left(
(-\Delta)^{-\frac{\alpha}{2}} \P ( {\vec{\mathfrak{f}}_0} ) 
\right)_i  
\right|^2
\neq 
\int 
\left|
\left(
(-\Delta)^{-\frac{\alpha}{2}} \P ( {\vec{\mathfrak{f}}_0} ) 
\right)_j
\right|^2$, or \\
\item[$ii)$] $\displaystyle  \int     \left(
(-\Delta)^{-\frac{\alpha}{2}} \P ( {\vec{\mathfrak{f}}_0} ) 
\right)_i  
\left(
(-\Delta)^{-\frac{\alpha}{2}} \P ( {\vec{\mathfrak{f}}_0} ) 
\right)_j
\neq 0$. 
\end{itemize}
In the following we will study both cases.

\medskip

{\bf The case $i)$}.  Note that 
\begin{equation*}
\begin{aligned}
\left|
\int \u_i^2 -  \u_j^2 - ( (\u_0)_i^2
+ (\u_0)_j^2)
\right| 
&
=
\left|\int\left(\u_i- (\u_0)_i\right) \u_i
+
\int \left(\u_i-(\u_0)_i\right)
(\u_0)_i
\right. 
\\
&
\quad \ \ 
\left.
+
\int\left(\u_j- (\u_0)_j\right) \u_j
+
\int \left(\u_j-(\u_0)_j\right)
(\u_0)_j
\right| 
\\
& \leq
C
\left\|
\u  - \u_0 \right\|_{L^2}
\left(
\|\u\|_{L^2}
+
\left\| \u_0\right\|_{L^2}
\right) .
\end{aligned}
\end{equation*}

Considering \eqref{sept20eq3} and \eqref{sept20eq2}
in the previous estimate, and recalling the notation $\u_0 = \eta  (-\Delta)^{-\frac{\alpha}{2}} \P ( {\vec{\mathfrak{f}}_0} )$, we get 
\begin{equation*}
    \begin{split}
 &\,  \left|
\int \u_i^2 -  \u_j^2 - 
\eta^2\left[\left(
(-\Delta)^{-\frac{\alpha}{2}} \P ( {\vec{\mathfrak{f}}_0} ) 
\right)_i^2 
+
\left(
(-\Delta)^{-\frac{\alpha}{2}} \P ( {\vec{\mathfrak{f}}_0} ) 
\right)_j^2
\right]
\right|  \\
\leq &\, C \eta^3
\left\|
(-\Delta)^{-\frac{\alpha}{2}} \P ( {\vec{\mathfrak{f}}_0} ) 
\right\|_{L^{ \frac{3}{\alpha-1}, \infty}}
\left\|
(-\Delta)^{-\frac{\alpha}{2}} \P ( {\vec{\mathfrak{f}}_0} ) 
\right\|_{L^2}^2 .
    \end{split}
\end{equation*}
Then, by imposing moreover 
\begin{equation*}
 \displaystyle   \eta_0
    \leq
    \frac{  
    \left| 
    \displaystyle 
\int
\left(
(-\Delta)^{-\frac{\alpha}{2}} \P ( {\vec{\mathfrak{f}}_0} ) 
\right)_i^2
+
\int
\left(
(-\Delta)^{-\frac{\alpha}{2}} \P ( {\vec{\mathfrak{f}}_0} ) 
\right)_j^2
\right|    
}{
\overline{C}  
\left\|
(-\Delta)^{-\frac{\alpha}{2}} \P ( {\vec{\mathfrak{f}}_0} ) 
\right\|_{L^{ \frac{3}{\alpha-1}, \infty}}
\left\|
(-\Delta)^{-\frac{\alpha}{2}} \P ( {\vec{\mathfrak{f}}_0} ) 
\right\|_{L^2}^2
}
,
\end{equation*}
with $\overline{C}>0$ selected  big enough, we conclude 
\begin{equation*}
    \begin{aligned}
        \left|
\int \u_i^2 -  \u_j^2 
\right|
& 
\geq
\eta^2 
    \left| 
    \displaystyle 
\int
\left(
(-\Delta)^{-\frac{\alpha}{2}} \P ( {\vec{\mathfrak{f}}_0} ) 
\right)_i^2
+
\int
\left(
(-\Delta)^{-\frac{\alpha}{2}} \P ( {\vec{\mathfrak{f}}_0} ) 
\right)_j^2
\right|     
\\
& \qquad -
\left|
\int \u_i^2 -  \u_j^2 - 
\eta^2 \left[
\int
\left(
(-\Delta)^{-\frac{\alpha}{2}} \P ( {\vec{\mathfrak{f}}_0} ) 
\right)_i^2 
+
\int
\left(
(-\Delta)^{-\frac{\alpha}{2}} \P ( {\vec{\mathfrak{f}}_0} ) 
\right)_j^2
\right]\right|
\\
& \geq 
\eta^2 
    \left| 
    \displaystyle 
\int
\left(
(-\Delta)^{-\frac{\alpha}{2}} \P ( {\vec{\mathfrak{f}}_0} ) 
\right)_i^2
+
\int
\left(
(-\Delta)^{-\frac{\alpha}{2}} \P ( {\vec{\mathfrak{f}}_0} ) 
\right)_j^2
\right|   
\\
& \qquad
-
12 C \eta^3
\left\|
(-\Delta)^{-\frac{\alpha}{2}} \P ( {\vec{\mathfrak{f}}_0} ) 
\right\|_{L^{ \frac{3}{\alpha-1}, \infty}}
\left\|
(-\Delta)^{-\frac{\alpha}{2}} \P ( {\vec{\mathfrak{f}}_0} ) 
\right\|_{L^2}^2
\\
& 
> 0.
    \end{aligned}
\end{equation*}
{\bf The case $ii)$}.  Note that 
\begin{equation*}
\begin{aligned}
\left|
\int \u_i \u_j-\int (\u_0)_i (\u_0)_j
\right| 
&
=
\left|\int\left(\u_i- (\u_0)_i\right) \u_j
+
\int (\u_0)_i\left(\u_j-(\u_0)_j\right)\right| 
\\
& \leq
\left\|
\u  - \u_0 \right\|_{L^2}
\left(
\|\u\|_{L^2}
+
\left\| \u_0\right\|_{L^2}
\right) .
\end{aligned}
\end{equation*}
Considering again \eqref{sept20eq3} and \eqref{sept20eq2}
in the expression above, we obtain 
\begin{equation*}
\begin{split}
&\,  \left|
\int \u_i \u_j-
\eta^2
\int
\left(
(-\Delta)^{-\frac{\alpha}{2}} \P ( {\vec{\mathfrak{f}}_0} ) 
\right)_i  
\left(
(-\Delta)^{-\frac{\alpha}{2}} \P ( {\vec{\mathfrak{f}}_0} ) 
\right)_j
\right|  \\
\leq &\, 12 C \eta^3
\left\|
(-\Delta)^{-\frac{\alpha}{2}} \P ( {\vec{\mathfrak{f}}_0} ) 
\right\|_{L^{ \frac{3}{\alpha-1}, \infty}}
\left\|
(-\Delta)^{-\frac{\alpha}{2}} \P ( {\vec{\mathfrak{f}}_0} ) 
\right\|_{L^2}^2.
\end{split}
\end{equation*}
Then, setting $\overline{C}>0$ big enough and by imposing moreover 
\begin{equation*}
 \displaystyle   \eta_0
    \leq
    \frac{  
    \left| 
    \displaystyle 
\int
\left(
(-\Delta)^{-\frac{\alpha}{2}} \P ( {\vec{\mathfrak{f}}_0} ) 
\right)_i  
\left(
(-\Delta)^{-\frac{\alpha}{2}} \P ( {\vec{\mathfrak{f}}_0} ) 
\right)_j
\right|    
}{
\overline{C}  
\left\|
(-\Delta)^{-\frac{\alpha}{2}} \P ( {\vec{\mathfrak{f}}_0} ) 
\right\|_{L^{ \frac{3}{\alpha-1}, \infty}}
\left\|
(-\Delta)^{-\frac{\alpha}{2}} \P ( {\vec{\mathfrak{f}}_0} ) 
\right\|_{L^2}^2
}
,
\end{equation*}
we can write
\begin{equation*}
    \begin{aligned}
        \left|
\int \u_i \u_j
\right|
& 
\geq
\eta^2 
    \left| 
    \displaystyle 
\int
\left(
(-\Delta)^{-\frac{\alpha}{2}} \P ( {\vec{\mathfrak{f}}_0} ) 
\right)_i  
\left(
(-\Delta)^{-\frac{\alpha}{2}} \P ( {\vec{\mathfrak{f}}_0} ) 
\right)_j
\right|    
\\
&
\qquad 
-
\left|
\int \u_i \u_j-
\eta^2
\int
\left(
(-\Delta)^{-\frac{\alpha}{2}} \P ( {\vec{\mathfrak{f}}_0} ) 
\right)_i  
\left(
(-\Delta)^{-\frac{\alpha}{2}} \P ( {\vec{\mathfrak{f}}_0} ) 
\right)_j
\right|
\\
& \geq 
\eta^2 
    \left| 
    \displaystyle 
\int
\left(
(-\Delta)^{-\frac{\alpha}{2}} \P ( {\vec{\mathfrak{f}}_0} ) 
\right)_i  
\left(
(-\Delta)^{-\frac{\alpha}{2}} \P ( {\vec{\mathfrak{f}}_0} ) 
\right)_j
\right|  \\
&
\qquad
- C \eta^3
\left\|
(-\Delta)^{-\frac{\alpha}{2}} \P ( {\vec{\mathfrak{f}}_0} ) 
\right\|_{L^{ \frac{3}{\alpha-1}, \infty}}
\left\|
(-\Delta)^{-\frac{\alpha}{2}} \P ( {\vec{\mathfrak{f}}_0} ) 
\right\|_{L^2}^2 
\\ 
&  
> 0.
    \end{aligned}
\end{equation*}
With this we conclude the proof of Proposition \ref{prop 10oct1}. \finpv

\medskip

In our last step, we deduce the follow estimates for this solution.
\begin{Proposition} Let $\u$ be the solution obtained in Proposition \ref{prop 10oct1}. Then, it holds: 
 \begin{equation}\label{strategy_lproof}
    C_1\left(  \frac{x}{|x|}  \right) 
    \frac{1}{|x|^{4-\alpha}}
    \leq 
|\u (x)|
    \leq 
    \frac{C_2}{|x|^{4-\alpha}}, \quad |x|\gg 1,
\end{equation} 
where $C_2>0$ and $C_1\left(  \frac{x}{|x|}  \right) >0$  is bounded for $|x| \gg 1$.   
\end{Proposition}
\pv  Let us start by proving the lower bound in (\ref{strategy_lproof}). We come back to the  asymptotic profile (\ref{sep 14_3}), where  for the sake of simplicity we shall denote 
\[ R_\alpha(x)= \begin{cases}
     O\left(
\frac{1}{|x|^{9-3\alpha}}
    \right) & \text{if }\alpha \neq 2,
     \\[9pt]
       O\left(
\frac{\log |x|}{|x|^{5-\alpha}}
    \right) & \text{if }\alpha = 2.
\end{cases}. \]
Then, for $|x|\gg 1$ we write 
\begin{equation}\label{Estim-below-01}
\left| m_\alpha(x):\int_{\Rt} (\u \otimes \u)(y)dy  \right| - \left| (-\Delta)^{-\frac{\alpha}{2}}\P(\f)(x)+ R_\alpha(x)\right| \leq |\u(x)|.    
\end{equation}
On the one hand, by Lemma \ref{Lem-Bilinear} the Kernel $m_\alpha(x)$ is homogeneous of degree $\alpha-4$ and we obtain
\begin{equation}\label{Estil-below-02}
    \begin{split}
   \left| m_\alpha(x):\int_{\Rt} (\u \otimes \u)(y)dy  \right|= &\,  \left| m_\alpha\left(|x| \frac{x}{|x|}\right):\int_{\Rt} (\u \otimes \u)(y)dy \right|\\
   = &\,  \frac{1}{|x|^{4-\alpha}} \,  \left| m_\alpha\left( \frac{x}{|x|}\right):\int_{\Rt} (\u \otimes \u)(y)dy \right| \\
   =&\, \frac{1}{|x|^{4-\alpha}} C_1 \left( \frac{x}{|x|} \right).
    \end{split}
\end{equation}
At this point, recall that always by Lemma \ref{Lem-Bilinear}  we have $\left| m_\alpha\left( \frac{x}{|x|}\right)\right|\leq C$, which yields
\[ C_1 \left( \frac{x}{|x|} \right) \leq C \, \left|\int_{\Rt} (\u \otimes \u)(y)dy \right| , \quad \mbox{for}\ \ |x|\gg 1.  \]
On the other hand, by our assumptions on the external force $\f$ given in Proposition \ref{prop 10oct1}, and by the expression $R_\alpha(x)$ defined above, we get
\begin{equation}\label{Estim-below-03}
\left| (-\Delta)^{-\frac{\alpha}{2}}\P(\f)(x)+ R_\alpha(x)\right|  = o \left( \frac{1}{|x|^{4-\alpha}}\right), \quad \mbox{for}\ \ |x|\gg 1.    
\end{equation}
With estimates (\ref{Estil-below-02}) and (\ref{Estim-below-03}) at hand, we get back to (\ref{Estim-below-01}) hence we obtain the lower bound in (\ref{strategy_lproof}).

\medskip

The upper bound in (\ref{strategy_lproof}) can be directly obtained by the estimate $|m_\alpha (x)|\leq C \frac{1}{|x|^{4-\alpha}}$ (see again Lemma \ref{Lem-Bilinear}) and by the good decaying propertied of  $\f$ and $R_\alpha$ already mentioned above. \finpv

\medskip

As already explained, the estimate from below in (\ref{strategy_lproof}) yields the statements of  Theorem \ref{Prop-non-existence}.

\section{Regularity of weak solutions: proof of Theorem \ref{Th3}}\label{section6}
For the sake of clearness, we shall divide the proof of Theorem \ref{Th3} into three main steps.

\medskip

{\bf Step 1. The parabolic setting.} Our starting point is to study the time-dependent fractional Navier-Stokes equations:
\begin{equation}\label{NS-Fractional-Parabolic}
\begin{cases}\vspace{2mm} 
\partial_t \v + (-\Delta)^{\frac{\alpha}{2}} \v +\P (\v \cdot \vec{\nabla})\v  = \P(\f), \quad   \text{div}(\v)=0, \quad  \alpha >1,\\
\v(0,\cdot)=\v_0,
\end{cases}
\end{equation}
where $\v_0$ denotes the (divergence-free) initial datum.  For a time $0<T<+\infty$, we denote  $\mathcal{C}_{*}([0,T],L^{p}(\Rt))$ the functional space of bounded and  weak$-*$ continuous  functions from $[0,T]$ with values in the  space $L^{p}(\Rt)$. Then, we shall prove the following:
\begin{Proposition}\label{Prop-NS-Parabolic} Let $1<\alpha$ and let $\max\left(\frac{3}{\alpha-1}, \frac{6}{\alpha},1\right)<p<+\infty$. Moreover, let   $\f\in L^p(\Rt)$ and $\v_0 \in  L^{p}(\Rt)$ be the external force and  the initial data respectively. There exists a time $T_0>0$,  depending on $\v_0$ and  $\f$, and there exists  a unique solution  $ \v  \in \mathcal{C}_{*}([0,T_0], L^{p}(\Rt))$ to equation (\ref{NS-Fractional-Parabolic}).  Moreover this solution verifies:   
	\begin{equation}\label{Cond-Sup}
	\sup_{0<t<T_0} t^{\frac{3}{\alpha  p}} \Vert  \v(t, \cdot) \Vert_{L^{\infty}} <+\infty. 
	\end{equation} 
\end{Proposition}
\noindent \pv

\noindent The proof is rather standard, so we shall detail the main estimates. Mild  solutions of the system (\ref{NS-Fractional-Parabolic}) write down as the integral formulation:
\begin{equation}\label{eqv}
\v(t,\cdot)= \,\, e^{-t (-\Delta)^{\frac{\alpha}{2}}}\v_0  + \int_{0}^{t} e^{-(t-s)(-\Delta)^{\frac{\alpha}{2}}} \,  \P(\f )ds  + \int_{0}^{t} e^{-(t-s)(-\Delta)^{\frac{\alpha}{2}}}\P (\text{div}(\v \otimes \v ))(s,\cdot) ds, 
\end{equation}
where we shall denote 
\begin{equation}\label{Bilinear-time}
\mathcal{B}(\v,\v)= \int_{0}^{t} e^{-(t-s)(-\Delta)^{\frac{\alpha}{2}}}\P (\text{div}(\v \otimes \v ))(s,\cdot) ds. 
\end{equation}
By the  Picard's fixed point argument,  we will solve the  problem  (\ref{eqv}) in the Banach space 
\[E_T = \left\{  g \in \mathcal{C}_{*}([0,T],L^{p}(\Rt)): \sup_{0<t<T} t^{\frac{3}{\alpha p}} \Vert g(t,\cdot)\Vert_{L^{\infty}}<+\infty \right\}, \]
with the norm 
\[ \Vert g \Vert_{E_T}= \sup_{0\leq t \leq T} \Vert g(t,\cdot)\Vert_{L^{p}}+\sup_{0<t<T} t^{\frac{3}{\alpha p}} \Vert g(t,\cdot)\Vert_{L^{\infty}}. \] 
We start by studying the  terms involving the data in  equation (\ref{eqv}). First, for the initial datum  $\v_0 \in L^{p}(\Rt)$   we have $\ds{\left\Vert e^{-t(-\Delta)^\frac{\alpha}{2}} \v_0 \right\Vert_{L^{p}} \leq c \Vert \v_0\Vert_{L^{p}}}$, hence we obtain  $\ds{e^{-t (-\Delta)^{\frac{\alpha}{2}}}\v_0 \in \mathcal{C}_{*}([0,T], L^p(\Rt))}$. Moreover,  by Lemma \ref{Lem-U2} we have $\ds{\sup_{0<t<T}t^{\frac{3}{\alpha p}} \left\Vert e^{-t(-\Delta)^\frac{\alpha}{2}} \v_0 \right\Vert_{L^{\infty}}  \leq c \Vert \v_0,  \Vert_{L^{p}}}$. We thus get $e^{-t (-\Delta)^\frac{\alpha}{2}}\v_0 \in E_T$ and it  holds:
\begin{equation}\label{Lin-Initial-Data}
\left\Vert e^{-t(-\Delta)^{\frac{\alpha}{2}}} \v_0  \right\Vert_{E_T} \leq c \Vert \v_0\Vert_{L^{p}}. 
\end{equation} 
Thereafter, for the external force $\f$ recall that  it is a  time-independent function. Then  we write
\begin{equation*}
\left\Vert \int_{0}^{t} e^{-(t-s)(-\Delta)^{\frac{\alpha}{2}}} \P (\f) ds  \right\Vert_{L^{p}}  \leq  \, \int_{0}^{t} \left\Vert e^{-(t-s)(-\Delta)^{\frac{\alpha}{2}}} \P(\f), \right\Vert_{L^{p}}\, ds 
\leq \,  c \left\Vert \f \right\Vert_{L^{p}} \left( \int_{0}^{t} ds\right),
\end{equation*} to get
\begin{equation*}\label{estim-forces-1}
\sup_{0 \leq t \leq T} \left\Vert \int_{0}^{t} e^{-(t-s)(-\Delta)^{\frac{\alpha}{2}}} \P(\f)  ds  \right\Vert_{L^{p}}  \leq c \, T\, \left\Vert  \f \right\Vert_{L^{p}}.
\end{equation*} 
On the other hand,   remark that by Lemma \ref{Lem-U2} we have
\[ \left\Vert e^{-(t-s)(-\Delta)^\frac{\alpha}{2}} \P(\f) \right\Vert_{L^{\infty}} \leq c\, (t-s)^{-\frac{3}{\alpha p}} \left\Vert \f  \right\Vert_{L^{p}},\]
and then we can  write  
\begin{equation*}
\begin{split}
&t^{\frac{3}{\alpha p}}\, \left\Vert \int_{0}^{t} e^{-(t-s)(-\Delta)^\frac{\alpha}{2}} \P(\f) ds \right\Vert_{L^\infty}  \leq t^{\frac{3}{\alpha p}} \, \int_{0}^{t} \left\Vert e^{-(t-s)(-\Delta)^\frac{\alpha}{2}}\P (\f) \right\Vert_{L^{\infty}} ds \\
 \leq &\,  c\, t^{\frac{3}{\alpha p}} \, \int_{0}^{t} (t-s)^{-\frac{3}{\alpha p}} \left\Vert \f  \right\Vert_{L^{p}} ds \leq c\,  t^{\frac{3}{\alpha p}}  \,  \left\Vert  \f\right\Vert_{L^{p}}  \, \left( \int_{0}^{t} (t-s)^{-\frac{3}{\alpha p}} \, ds  \right) \leq c \, t\, \left\Vert \f  \right\Vert_{L^{p}}. 
\end{split}
\end{equation*}
We thus obtain
\begin{equation*}\label{estim-forces-2}
\sup_{0 <t < T} t^{\frac{3}{\alpha p}} \, \left\Vert \int_{0}^{t} e^{-(t-s)(-\Delta)^\frac{\alpha}{2}}\P(\f) ds \right\Vert_{L^\infty}  \leq c \, T\, \left\Vert \f  \right\Vert_{L^{p}}.
\end{equation*}
By the estimates above  we get
\begin{equation}\label{estim-forces}
\left\Vert \int_{0}^{t} e^{-(t-s)(-\Delta)^\frac{\alpha}{2}} \P(\f)ds \right\Vert_{E_T}  \leq c\, T\, \left\Vert  \f  \right\Vert_{L^{p}}.
\end{equation}
Now, we study the bilinear form $\mathcal{B}(\v,\v)$ defined in (\ref{Bilinear-time}). Our starting point is to prove the estimate 
\begin{equation}\label{Estim-Bilin1}
\sup_{0\leq t \leq T} \left\Vert \mathcal{B}(\v, \v)  \right\Vert_{L^{p}}  \leq  c\, T^{1-\frac{1}{\alpha}-\frac{3}{\alpha p}}\, \Vert \v \Vert^{2}_{E_T}, \qquad  1-\frac{1}{\alpha}-\frac{3}{\alpha p}>0,
\end{equation} 
where remark that $1-\frac{1}{\alpha}-\frac{3}{\alpha p}>0$ as long as $p>\frac{3}{\alpha-1}$. Indeed, we have
\begin{equation*}\label{BiLin1}
\begin{split}
 \sup_{0\leq t \leq T} \left\Vert \mathcal{B}(\v, \v)  \right\Vert_{L^{p}} \leq &\,\,  c\,  \sup_{0\leq t \leq T} \int_{0}^{t}  \left\Vert   e^{-(t-s)(-\Delta)^{\frac{\alpha}{2}}} (\text{div}(\v \otimes \v ))(s,\cdot)\right\Vert_{L^{p}}  ds\\
\leq & \,\,c\,  \sup_{0\leq t \leq T} \int_{0}^{t} \frac{1}{(t-s)^{\frac{1}{\alpha}}} \Vert \v(s,\cdot) \otimes  \v(s,\cdot)\Vert_{L^{p}}\, ds\\ 
\leq &\,\, c\, \sup_{0\leq t \leq T} \int_{0}^{t} \frac{1}{(t-s)^{\frac{1}{\alpha}}\,s^{\frac{3}{\alpha p}} }   (s^{\frac{3}{\alpha p}}\Vert \v(s,\cdot) \Vert_{L^{\infty}}) \Vert \v(s,\cdot)\Vert_{L^{p}}\,ds\\
\leq & \,\,c\, T^{1-\frac{1}{\alpha}-\frac{3}{\alpha p}}\,  \left\Vert  \v  \right\Vert^{2}_{E_T}. 
\end{split}
\end{equation*}  
Thereafter, we shall prove the estimate:
\begin{equation}\label{BiLin2}
\sup_{0\leq t \leq T} t^{\frac{3}{\alpha p}} \left\Vert \mathcal{B}(\v, \v)  \right\Vert_{L^{\infty}} \leq  c\, T^{1-\frac{1}{\alpha}-\frac{3}{\alpha p}}\,  \left\Vert \v  \right\Vert^{2}_{E_T}, \qquad 1-\frac{1}{\alpha}-\frac{3}{\alpha p}>0. 
\end{equation}
We write
\begin{equation*}
\begin{split}
 \sup_{0\leq t \leq T} t^{\frac{3}{\alpha p}} \left\Vert \mathcal{B}(\v, \v)  \right\Vert_{L^{\infty}}  \leq  \,\,  \sup_{0<t<T} t^{\frac{3}{\alpha p}}    \int_{0}^{t}  \left\Vert  e^{-(t-s)(-\Delta)^{\frac{\alpha}{2}}} \P \, \text{div} \left( \v \otimes \v  \right)(s,\cdot)  \right\Vert_{L^{\infty}} ds=(a).
\end{split}
\end{equation*}
Then, by the third point of Lemma \ref{Lem-Kernel}, and using now the information $p>\frac{6}{\alpha}$,  we can write:
\begin{equation*}
\begin{split}
(a)\leq &  \,\, c\,  \sup_{0\leq t \leq T} t^{\frac{3}{\alpha p}}   \int_{0}^{t} \frac{1}{(t-s)^{\frac{1}{\alpha}}}  \Vert \v(s,\cdot) \otimes  \v(s,\cdot)\Vert_{L^{\infty}} ds\\
\leq & \,\,  c\,  \sup_{0\leq t \leq T} t^{\frac{3}{\alpha p}}   \int_{0}^{t} \frac{ds}{(t-s)^{\frac{1}{\alpha}} s^{\frac{6}{\alpha p}}} \left( s^{\frac{3}{\alpha p}} \Vert \v(s,\cdot) \Vert_{L^{\infty}}\right)^2 ds\\
\leq & \,\,  c\, \left( \sup_{0\leq t \leq T} t^{\frac{3}{\alpha p}}    \int_{0}^{t} \frac{ds}{(t-s)^{\frac{1}{\alpha}} s^{\frac{6}{\alpha p}}} \right) \left\Vert  \v  \right\Vert^{2}_{E_T}. \\
\leq & \,\,  c \left( \sup_{0\leq t \leq T}\left[ t^{\frac{3}{\alpha p}}   \int_{0}^{t/2} \frac{ds}{(t-s)^{\frac{1}{\alpha}} s^{\frac{6}{\alpha p}}} +t^{\frac{3}{\alpha p}}   \int_{t/2}^{t} \frac{ds}{(t-s)^{\frac{1}{\alpha}} s^{\frac{6}{\alpha p}}} \right]\right)  \left\Vert \v  \right\Vert^{2}_{E_T}\\
\leq &  \,\,  c\,\left(  \sup_{0\leq t \leq T}\left[ t^{\frac{3}{\alpha p}-\frac{1}{\alpha}}    \int_{0}^{t/2} \frac{ds}{s^{\frac{6}{\alpha p}}} +t^{\frac{3}{\alpha p}-\frac{6}{\alpha p}}  \int_{t/2}^{t} \frac{ds}{(t-s)^{\frac{1}{\alpha}}} \right]\right)  \left\Vert \v  \right\Vert^{2}_{E_T}\\
\leq &  \,\,  c\, T^{1-\frac{1}{\alpha}-\frac{3}{\alpha p}}\,  \left\Vert \v  \right\Vert^{2}_{E_T}.
\end{split} 
\end{equation*} 
By the inequalities (\ref{BiLin1}) and (\ref{BiLin2}) we can write:
\begin{equation}\label{Bilinear}
\| \mathcal{B}(\v, \v) \|_{E_T} \leq  c\, T^{1-\frac{1}{\alpha}-\frac{3}{\alpha p}}\, \| \v\|^{2}_{E_T}, \qquad 1-\frac{1}{\alpha}-\frac{3}{\alpha p}>0.    
\end{equation}

Once we have the estimates (\ref{Lin-Initial-Data}), (\ref{estim-forces}) and (\ref{Bilinear}) at our disposal, the proof of Proposition \ref{Prop-NS-Parabolic}  follows from well-known arguments. \finpv 

\medskip 

{\bf Step 2. Global boundedness of $\u$}. With the  help of the Proposition \ref{Prop-NS-Parabolic},  we are able to prove  the following:  
\begin{Proposition}\label{Prop-Global-boundedness} Let $1<\alpha$, $\max\left(\frac{3}{\alpha-1},\frac{6}{\alpha},1 \right)<p<+\infty$ and $0\leq s$. Let $\f \in \dot{W}^{-1,p}\cap \dot{W}^{s,p}(\Rt)$ be the external force and let   $\u \in  L^{p}(\Rt)$  be a weak solution to equation (\ref{Frac-NS}). Then we have $\u \in L^{\infty}(\Rt)$.  
\end{Proposition} 	
\noindent \pv

\noindent First remark that by hypothesis on the external force and by interpolation (see \cite[Chapter $5$]{adams2003pseudo}) we have $\f \in L^p(\Rt)$. 

\medskip 

Then, in the initial value  problem (\ref{NS-Fractional-Parabolic}), we set  the initial data $\ds{\v_0=\u}$. Then, by  Proposition \ref{Prop-NS-Parabolic} there exists a time $0<T_0$, and there exists  a unique  solution $\ds{\v  \in \mathcal{C}_{*}([0,T_0], L^{p}(\Rt))}$ to  equation  (\ref{NS-Fractional-Parabolic}), which arises from $\u$.

\medskip 

On the other hand, we have the following key remark: since $\u$ is a  time-independent function we have $\partial_t \u =0$ and this function is also a solution of the initial value problem (\ref{NS-Fractional-Parabolic}) with the same  initial data $\u$. Moreover, we have $\ds{\u \in \mathcal{C}_{*}([0,T_0], L^{p}(\Rt))}$.

\medskip  

Consequently, in the space $\ds{\mathcal{C}_{*}([0,T_0], L^{p}(\Rt))}$ we  have   two solutions of  equation  (\ref{NS-Fractional-Parabolic}) with the same initial data $\v_0=\u$: the solution $\v$ given by the Proposition \ref{Prop-NS-Parabolic} and the time-independent solution $\u$. By  uniqueness we have the identity $\v=\u$  and by   (\ref{Cond-Sup})  we can write 
\[  \sup_{0<t<T_0} t^{\frac{3}{2 p}} \, \Vert  \u \Vert_{L^{\infty}} <+\infty. \]
But, as the solution $\u$  does not depend on the temporal variable we finally get $\u \in L^{\infty}(\Rt)$ and Proposition  \ref{Prop-Global-boundedness} is now proven. \finpv  

\medskip

{\bf Step 3. Regularity of $\u$ and $P$.}  The global boundedness of $\u$ obtained in the last step is the key tool to prove the following: 
\begin{Proposition}\label{Prop-Regularity} Since $\f \in  \dot{W}^{-\alpha,p}\cap \dot{W}^{s,p}(\Rt)$, with $1<\alpha$, $0\leq s$ and $\max\left(\frac{3}{\alpha-1},\frac{6}{\alpha},1\right)<p<+\infty$, and since $\u \in L^p \cap L^{\infty}(\Rt)$, we have $\u \in \dot{W}^{s+\alpha,p}(\Rt)$ and $P \in \dot{W}^{s+\alpha,p}(\Rt) +  \dot{W}^{s+1,p}(\Rt)$. 
\end{Proposition}
\noindent \pv

\noindent Let us start by explaining the general strategy of the proof.  For $0\leq s$ and $1<\alpha$, we consider the quantities $0< s+\alpha$ and $0<\alpha-1$. Then, let $k \in \mathbb{N}$ be such that $k(\alpha-1) \leq s+\alpha \leq (k+1)(\alpha-1)$. We thus  write $s+\alpha= k(\alpha-1)+\varepsilon$, with $0 \leq \varepsilon < \alpha-1$. To prove that $\u \in \dot{W}^{s+\alpha,p}(\Rt)$, first we shall prove that $\u \in \dot{W}^{k(\alpha-1),p}(\Rt)$  and next we will verify that $(-\Delta)^{\frac{k(\alpha-1)}{2}}\u \in \dot{W}^{\varepsilon,p}(\Rt)$.  

\medskip

By an iteration process, we prove that $\u \in \dot{W}^{k(\alpha-1),p}(\Rt)$. For the sake of simplicity, we consider the following cases of $k$. 

\medskip

First, if  $k=0$ by our hypothesis we directly have $\u \in L^p(\Rt)$. Next, if $k=1$, recall that $\u$ solves the fixed point equation (\ref{Frac-NS-fixed-point}). Then we have 
\begin{equation*}
(-\Delta)^{\frac{\alpha-1}{2}} \u = (-\Delta)^{-\frac{1}{2}} \text{div}(\u \otimes \u) +(-\Delta)^{-\frac{1}{2}} \P (\f).     
\end{equation*}
For the first term on the right-hand side, recall that $\u \in L^p\cap L^\infty(\Rt)$ and we get $(-\Delta)^{-\frac{1}{2}} \text{div}(\u \otimes \u) \in L^p(\Rt)$. For the second term on the right-hand side, by our hypothesis $\f \in \dot{W}^{-\alpha,p}\cap \dot{W}^{s,p}(\Rt)$ with $1<\alpha$ and $0\leq s$ and  by interpolation (see \cite[Chapter $5$]{adams2003pseudo})  we have 
 $(-\Delta)^{-\frac{1}{2}} \P (\f) \in L^p(\Rt)$. We thus get $\u \in \dot{W}^{\alpha-1,p}(\Rt)$.   
 
 \medskip

Thereafter, if $2\leq k$, we start by  applying the operator $(-\Delta)^{\frac{\alpha-1}{2}}$ to the identity and we get 
\begin{equation*}
(-\Delta)^{\frac{2(\alpha-1)}{2}} \u = (-\Delta)^{\frac{\alpha-1}{2}} (-\Delta)^{-\frac{1}{2}} \text{div}(\u \otimes \u) + (-\Delta)^{\frac{\alpha-1}{2}} (-\Delta)^{-\frac{1}{2}} \P (\f).     
\end{equation*}
As before, we shall prove that each term on the right-side belong to the space $L^p(\Rt)$. For the first term, we write $(-\Delta)^{\frac{\alpha-1}{2}} (-\Delta)^{-\frac{1}{2}} \text{div}(\u \otimes \u)= (-\Delta)^{-\frac{1}{2}} \text{div}((-\Delta)^{\frac{\alpha-1}{2}} ( \u \otimes \u))$. Thereafter, since $(-\Delta)^{\frac{\alpha-1}{2}} \u \in L^p(\Rt)$ and $\u \in L^\infty(\Rt)$ by Lemma \ref{Leibniz-rule} we write $\| (-\Delta)^{\frac{\alpha-1}{2}} (\u \otimes \u) \|_{L^p} \leq c \| (-\Delta)^{\frac{\alpha-1}{2}} \u \|_{L^p} \| \u \|_{L^\infty}$, and we have $(-\Delta)^{\frac{\alpha-1}{2}} (-\Delta)^{-\frac{1}{2}} \text{div}(\u \otimes \u) \in L^p(\Rt)$.  For the second term we write  $(-\Delta)^{\frac{\alpha-1}{2}} (-\Delta)^{-\frac{1}{2}} \P (\f)=(-\Delta)^{\frac{2(\alpha-1)-\alpha}{2}} \P (\f)$. Then, as $2(\alpha-1)-\alpha \leq k(\alpha-1)-\alpha \leq s$ (recall that $s+\alpha=k(\alpha-1)+\varepsilon$), and moreover, as $\f \in \dot{W}^{-\alpha,p}\cap \dot{W}^{s,p}(\Rt)$,  we have $(-\Delta)^{\frac{2(\alpha-1)-\alpha}{2}} \P (\f) \in L^p(\Rt)$. We thus get $\u \in \dot{W}^{2(\alpha-1),p}(\Rt)$.  By iterating this process until $k$ we obtain that $\u \in \dot{W}^{k(\alpha-1),p}(\Rt)$.

\medskip

Finally, to prove that $\u \in \dot{W}^{s+\alpha,p}(\Rt)$ we must verify that $(-\Delta)^{\frac{\varepsilon}{2}} (-\Delta)^{\frac{k(\alpha-1)}{2}} \u \in L^p(\Rt)$. We thus write 
\begin{equation*}
\begin{split}
(-\Delta)^{\frac{\varepsilon}{2}} (-\Delta)^{\frac{k(\alpha-1)}{2}} \u = &\,  (-\Delta)^{\frac{\varepsilon}{2}} \left( (-\Delta)^{\frac{(k-1)(\alpha-1)}{2}} (-\Delta)^{-\frac{1}{2}} \text{div}(\u \otimes \u) + (-\Delta)^{\frac{(k-1)(\alpha-1)}{2}}(-\Delta)^{-\frac{1}{2}} \P(\f)\right)\\
=&\, (-\Delta)^{\frac{\varepsilon+(k-1)(\alpha-1)}{2}} (-\Delta)^{-\frac{1}{2}}\text{div}(\u \otimes \u) + (-\Delta)^{\frac{\varepsilon+(k-1)(\alpha-1) -1}{2}}\P(\f). 
\end{split}
\end{equation*}
For the first term on the right-hand side, since $(-\Delta)^{\frac{k(\alpha-1)}{2}} \u \in L^p(\Rt)$ and $\u \in L^\infty(\Rt)$, by Lemma \ref{Leibniz-rule} we have  $(-\Delta)^{\frac{k(\alpha-1)}{2}} (-\Delta)^{-\frac{1}{2}}\text{div}(\u \otimes \u) \in L^p(\Rt)$.  Moreover, since $(-\Delta)^{\frac{(k-1)(\alpha-1)}{2}} (-\Delta)^{-\frac{1}{2}}\text{div}(\u \otimes \u) \in L^p(\Rt)$ and $0 \leq  \varepsilon < \alpha-1$, by interpolation we obtain  $(-\Delta)^{\frac{\varepsilon+(k-1)(\alpha-1)}{2}} (-\Delta)^{-\frac{1}{2}}\text{div}(\u \otimes \u)\in L^p(\Rt)$. On the other hand, for the second term on the right-hand side, recall that $s+\alpha= k(\alpha-1)+\varepsilon$, hence $\varepsilon+(k-1)(\alpha-1) -1=s$. Then, by our hypothesis $\f \in \dot{W}^{s,p}(\Rt)$ we  directly  have $(-\Delta)^{\frac{\varepsilon+(k-1)(\alpha-1) -1}{2}}\P(\f) \in L^p(\Rt)$. We thus obtain $\u \in \dot{W}^{s+\alpha,p}(\Rt)$.
\begin{Remarque}\label{Rmk-regularity} The initial regularity $\f \in \dot{W}^{s,p}(\Rt)$ stops this iterative process  yielding  that $\u \in \dot{W}^{s+\alpha,p}(\Rt)$ is the \emph{maximum} gain of regularity for solutions.  
\end{Remarque}

Now, we study the pressure term $P$. Recall that $P$ is related to the velocity $\u$ and the external force $\f$ by the expression 
\begin{equation*}
 P= (-\Delta)^{-1} \text{div}(\text{div}(\u \otimes \u)) - (-\Delta)^{-1}\text{div}(\f).   
\end{equation*}
For the first term on the right-hand side, since $\u \in \dot{W}^{s+\alpha,p}(\Rt)$ and $\u \in L^\infty(\Rt)$, by Lemma \ref{Leibniz-rule} we get $(-\Delta)^{-1} \text{div}(\text{div}(\u \otimes \u)) \in \dot{W}^{s+\alpha,p}(\Rt)$. Moreover, for the second term on the right-hand side,  we directly have $(-\Delta)^{-1}\text{div}(\f) \in \dot{W}^{s+1,p}(\Rt)$.  Proposition \ref{Prop-Regularity} is proven. \finpv

\medskip

With this we conclude the proof of
Theorem \ref{Th3}. \finpv

\section{Liouville-type result: Proof of Proposition \ref{Proposition-Liouville}}\label{section7}
The proof follows some of the  main estimates performed in \cite{chamorropoggi2023p} and \cite{chamorrojarrinlem2018p}. First, we consider $\varphi \in \mathcal{C}^{\infty}_{0}(\Rt)$ a positive and radial function such that $\varphi(x)=1$ when $|x|<1/2$ and $\varphi(x)=0$ when $|x|\geq 1$. Then, for $R\geq 1$, we define the cut-off function $\varphi_R(x)= \varphi(x/R)$. Remark that $supp(\varphi_R) \subset B_R$, where we denote $B_R=\{ x \in \Rt: \ |x| \leq R\}$. 

\medskip

On the other hand, since $\f \equiv 0$ and since $\u \in L^p(\Rt)$ with $\max\left(\frac{3}{\alpha-1},\frac{6}{\alpha}\right)<p$, by Theorem \ref{Th3} we have $\u \in \mathcal{C}^{\infty}(\Rt)$ and $P \in \mathcal{C}^{\infty}(\Rt)$. Therefore, we can multiply each term in  equation (\ref{Frac-NS-homog}) by $\varphi_R \u$, then  we integrate by parts over $\Rt$ to obtain the estimate proven in \cite[Estimate $(4.5)$]{chamorropoggi2023p}: 
\begin{equation}\label{Estim-local-liouville}
\begin{split}
\int_{B_{\frac{R}{2}}}| (-\Delta)^{\frac{\alpha}{4}} \u |^2 dx \leq &\,  \int_{B_R} \vec{\nabla}\varphi_R \cdot \left( \frac{|\u|^2}{2}+P\right)\u \,  dx \\
&\, + \int_{\Rt}  (-\Delta)^{\frac{\alpha}{4}} \u \cdot \left( \varphi_R ((-\Delta)^{\frac{\alpha}{4}} \u) -(-\Delta)^{\frac{\alpha}{4}}(\varphi_R \u)\right) dx\\
=&\, I_1+I_2.
\end{split}
\end{equation}
By  H\"older inequalities and the fact that $supp (\vec{\nabla}\varphi_R)\subset C(\frac R 2, R)$, where we denote $C(\frac R 2, R)=\{ x \in \Rt : \ \frac R 2 < |x|<R\}$, the term $I_1$ was estimated in \cite[Proof of Theorem $1$]{chamorrojarrinlem2018p} and we have  
\begin{equation*}
  I_1 \leq C R^{2-\frac{9}{p}} \| \vec{\nabla}\varphi \|_{L^r} \| \u \|^{3}_{L^p(C(\frac{R}{2},R))}, \quad 1=\frac{1}{r}+\frac{3}{p},   
\end{equation*}
and since $p\leq \frac{9}{2}$ we obtain 
\begin{equation*}
    I_1 \leq C\, \| \u \|^{3}_{L^p(C(\frac{R}{2},R))}. 
\end{equation*}
On the other hand, for $\frac{5}{3}<\alpha_1<\alpha$ and $0<\alpha_2<\alpha$ such that $\alpha_1+\alpha_2=\alpha\leq 2$, and for $1<p_1<+\infty$ such that $1/2=1/p_1 + 1/p$, the term $I_2$ was estimated in \cite[Page $13$]{chamorropoggi2023p} as follows:
\begin{equation*}
   I_2 \leq C \| \u \|_{\dot{H}^{\frac{\alpha}{2}}} \left( \| (-\Delta)^{\frac{\alpha_1}{4}} \varphi_R \|_{L^{p_1}}\, \| (-\Delta)^{\frac{\alpha_2}{4}} \u \|_{L^p}+\|(-\Delta)^{\frac{\alpha}{4}}\varphi_R \|_{L^{p_1}}\,\| \u \|_{L^p}\right).  
\end{equation*}
By the localization properties of $\varphi_R$ we can write 
\begin{equation*}
 I_2 \leq C \| \u \|_{\dot{H}^\frac{\alpha}{2}} \left(R^{-\frac{\alpha_1}{2}+\frac{3}{p_1}}\| (-\Delta)^{\frac{\alpha_2}{2}} \u \|_{L^p} + R^{-\frac{\alpha}{2}+\frac{3}{p_1}} \| \u \|_{L^p} \right).  
\end{equation*}
Then, since $\alpha_1<\alpha$ we have 
\begin{equation*}
   I_2 \leq C  \| \u \|_{\dot{H}^\frac{\alpha}{2}} \, R^{\frac{-\alpha_1}{2}+\frac{3}{p_1}} \left(\| (-\Delta)^{\frac{\alpha_2}{2}} \u \|_{L^p} + \| \u \|_{L^p} \right).
\end{equation*}
Here, remark that by hypothesis we have $\| \u \|_{L^p} <+\infty$, and since $\f \equiv 0$ by Theorem \ref{Th-Picard} we also have $\| (-\Delta)^{\frac{\alpha_2}{2}} \u \|_{L^p} <+\infty$. 
Gathering the estimates on $I_1$ and $I_2$, we get back to (\ref{Estim-local-liouville}) to obtain the following local estimate
\begin{equation*}
    \int_{B_{\frac{R}{2}}}|(-\Delta)^{\frac{\alpha}{4}}\u |^2 dx \leq C\, \| \u \|^{3}_{L^p(C(\frac{R}{2},R))}+ C  \| \u \|_{\dot{H}^\frac{\alpha}{2}} \, R^{\frac{-\alpha_1}{2}+\frac{3}{p_1}} \left(\| (-\Delta)^{\frac{\alpha_2}{2}} \u \|_{L^p} + \| \u \|_{L^p} \right).
\end{equation*}
For the fist term on the right-hand side, since $\u \in L^p(\Rt)$ we have $\ds{\lim_{R\to +\infty} \| \u \|^{3}_{L^p(C(\frac{R}{2},R))}=0}$.  For the second term on the right-hand side, we shall verify that $\frac{-\alpha_1}{2}+\frac{3}{p_1}<0$. Indeed, by the relationship $\frac 1 2=\frac 1 p_1+ \frac 1 p$ a simple computation shows that the inequality $\frac{-\alpha_1}{2}+\frac{3}{p_1}<0$ is equivalent to $p<\frac{6}{3-\alpha_1}$. Moreover, since $p\leq \frac{9}{2}$ this last inequality holds as long as $\frac 9 2<\frac{6}{3-\alpha_1}$, which is ultimately verified by the constrain $\frac 5 3<\alpha_1$. 
With this we conclude the proof of Proposition \ref{Proposition-Liouville}.\finpv 

\section*{Appendix: New regularity criterion for classical stationary Navier-Stokes equations in Morrey spaces}\label{Appendix}

Let consider the classical stationary Navier-Stokes equations in the whole space $\Rt$:
\begin{equation}\label{Classical-NS-2}
-\Delta  \u +  (\u \cdot \vec{\nabla}) \u + \vec{\nabla}P = \text{div}(\mathbb{F}), \quad \text{div}(\u)=0,
\end{equation}
where, for the sake of simplicity, we  have written the external force $\f$ as the divergence of a tensor  $\mathbb{F}=(F_{ij})_{1\leq i,j\leq 3}$.  

\medskip

We shall name a \emph{very weak} of equation (\ref{Classical-NS-2}) the couple $(\u,P)\in L^{2}_{loc}(\Rt)\times \mathcal{D}'(\Rt)$ which verifies this equation in the distributional sense.   Note that very weak solutions  have  minimal conditions to guarantee  that each term in equation (\ref{Classical-NS-2}) is well-defined as distributions.  In particular, we let the pressure $P$ to be a very general object as we
only have $P \in \mathcal{D}'(\Rt)$.

\medskip

As the velocity $\u$ is a locally square integrable function, in order to improve their regularity we look for some natural conditions on the local quantities $\ds{\int_{B(x_0,R)} |\u(x)|^2 dx}$, where $B(x_0,\R)=\{ x\in \Rt : \ |x-x_0|<R\}$. Thus, the Morrey spaces appear naturally. 

\medskip

Recall that  for $2<p<+\infty$ the homogeneous Morrey space $\dot{M}^{2,p}(\Rt)$ is the space of $L^{2}_{loc}$-functions such that 
\begin{equation*}
    \| f \|_{\dot{M}^{2,p}}= \sup_{R>0, \ x_0 \in \Rt} R^{\frac{3}{p}} \left( \frac{1}{dx(B(x_0,R))} \int_{B(x_0, R)} | f(x)|^2 dx  \right)^{\frac{1}{2}}<+\infty,
\end{equation*}
where $dx(B(x_0,R)) \sim R^3$ is the Lebesgue measure of the ball $B(x_0,R)$. Here, the parameter $p$ measures the decay rate of the local quantity $\ds{\left( \frac{1}{dx(B(x_0,R))} \int_{B(x_0, R)} | f(x)|^2 dx  \right)^{\frac{1}{2}}}$ when $R \to +\infty$. Moreover, this is a homogeneous space of order $-\frac{3}{p}$, and the following chain of continuous embeddings holds $L^p(\Rt)\subset L^{p,q}(\Rt)\subset \dot{M}^{2,p}(\Rt)$. We thus study the regularity of very weak solutions in a general framework.  

\medskip

For the parameter $p$ given above and for a regularity parameter $k\in \mathbb{N}$, we introduce the following Sobolev-Morrey space
\begin{equation*}
   \mathcal{W}^{k,p}=\big\{ f \in \dot{M}^{2,p}(\Rt): \ \partial^{{\bf a}} f \in \dot{M}^{2,p}(\Rt)\ \  \mbox{for all multi-indice} \ |{\bf a}|\leq k \big\}.  
\end{equation*}
Moreover, we denote by $W^{k,\infty}(\Rt)$ the classical Sobolev space of bounded functions with bounded weak derivatives until the other $k$.  Finally, for $0 < s < 1$, we shall denote by $\mathcal{C}^{k,s}(\Rt)$  the H\"older space of $\mathcal{C}^k$-
functions whose  derivatives are H\"older continuous functions with parameter $s$. 

\medskip 

In this framework, we obtain a new regularity criterion for \emph{very weak} solutions to equation (\ref{Classical-NS-2}).
\begin{Theoreme}\label{Th-Reg-Classical}
Let $(\u,P)\in L^{2}_{loc}(\Rt)\times \mathcal{D}'(\Rt)$ be a very weak solution to (\ref{Classical-NS-2}), with $\mathbb{F} \in \mathcal{D}'(\Rt)$.  For $k \geq  0$ and $3<p$ assume that 
\begin{equation}\label{Hyp-Force}
\mathbb{F}   \in \mathcal{W}^{k+1,p}(\Rt) \cap W^{k+1,\infty}(\Rt).
\end{equation}
Then, if the velocity field verifies 
\begin{equation*}
\u \in \dot{M}^{2,p}(\Rt),
\end{equation*}
it follows that $\u \in \mathcal{W}^{k+2,p}(\Rt)$ and $P \in \mathcal{W}^{k+1,p}(\Rt)$.  Moreover,  we have $\u \in \mathcal{C}^{k+1,s}(\Rt)$ and $P\in \mathcal{C}^{k,s}(\Rt)$ with $s=1-\frac{3}{p}$.    
\end{Theoreme} 

Recall that the external force acting on equation (\ref{Classical-NS-2}) is given by $\text{div}(\mathbb{F})$. Then, by our assumption (\ref{Hyp-Force}) we have $\text{div}(\mathbb{F}) \in \mathcal{W}^{k,p}$, which yields a gain of regularity of very weak solution of the order $k+2$. As in Theorem \ref{Th2}, this \emph{maximum} gain of regularity is given by the effects of the Laplacian operator in equation (\ref{Classical-NS-2}). 

\medskip

On the other hand, in the particular homogeneous case $\mathbb{F} \equiv 0$, very weak solutions to equation (\ref{Classical-NS-2}) verify $(\u, P) \in \mathcal{C}^{\infty}(\Rt)$, provided that $\u \in \dot{M}^{2,p}(\Rt)$ with $3<p$. As explained before introducing Proposition \ref{Proposition-Liouville}, this particular result is of interest in connection to the Liouville-type problem for equation (\ref{Classical-NS-2})  in the  setting of Morrey spaces \cite{chamorrojarrinlem2018p}.

\medskip

\noindent{\bf Proof of Theorem \ref{Th-Reg-Classical}}.\\
\noindent The proof follows the main ideas in the proof of Theorem \ref{Th3} (see Section  \ref{section6}), so we shall only detail the main computations. 

\medskip 

By well-known properties of Morrey spaces,  Proposition \ref{Prop-NS-Parabolic} also holds in the (larger) space $\dot{M}^{2,p}(\Rt)$  and we can state
\begin{Proposition} For $3<p$,  let   $\text{div}(\mathbb{F})\in \dot{M}^{2,p}(\Rt)$ and $\v_0 \in  \dot{M}^{2,p}(\Rt)$ be the external force and  the (divergence-free) initial data respectively. There exists a time $T_0>0$,  depending on $\v_0$ and  $\text{div}(\mathbb{F})$, and there exists  a unique solution $\v  \in \mathcal{C}_{*}([0,T_0], \dot{M}^{2,p}(\Rt))$ to the Navier-Stokes equations:
\begin{equation*}\
\partial_t \v  -\Delta \v +\P (\v \cdot \vec{\nabla})\v  = \P(\f), \quad   \text{div}(\v)=0, \quad 
\v(0,\cdot)=\v_0.
\end{equation*}
Moreover this solution verifies $\ds{	\sup_{0<t<T_0} t^{\frac{3}{2 p}} \Vert  \v(t, \cdot) \Vert_{L^{\infty}} <+\infty}$.
\end{Proposition}

Moreover, by following the same ideas in the proof of Proposition \ref{Prop-Global-boundedness}, we obtain our key result
\begin{Proposition}\label{Prop-Global-boundedness2} Let $3<p$ and $0\leq k$. Let $\mathbb{F} \in \mathcal{W}^{k+1,p}(\Rt)$ be the tensor  let   $\u \in  \dot{M}^{2,p}(\Rt)$  be a very weak solution to equation (\ref{Classical-NS-2}). Then we have $\u \in L^{\infty}(\Rt)$.  
\end{Proposition}

Finally,  global boundedness of the velocity $\u$ allows us to study its regularity. 
\begin{Proposition}\label{Prop-Regularity-Classical}
Since $\mathbb{F} \in  \mathcal{W}^{k+1,p}\cap W^{k,\infty}(\Rt)$, with $3<p$,  and since $\u \in \dot{M}^{2,p} \cap L^{\infty}(\Rt)$, we have $\u \in \mathcal{W}^{k+2,p}(\Rt)$ and $P \in \mathcal{W}^{k+1,p}(\Rt)$. Moreover, it holds  $\u \in \mathcal{C}^{k+1,s}(\Rt)$ and $P\in \mathcal{C}^{k,s}(\Rt)$ with $s=1-\frac{3}{p}$.  
\end{Proposition}
\noindent\pv\\
As in the proof of Proposition \ref{Prop-Regularity}, we consider the fixed point equation 
\begin{equation}\label{Fixed-point-NS-Classical}
    \u = -(-\Delta)^{-1} \P\big(\text{div}(\u \otimes \u)\big) + (-\Delta)^{-1}\P  (\text{div}(\mathbb{F})). 
\end{equation} 
By using this equation, we shall prove that $\partial^{{\bf a}} \u \in \dot{M}^{2,p}(\Rt)$ for all multi-indice $|{\bf a}|\leq k+2$.  We shall prove this fact by iteration respect to the order of the multi-indices ${\bf a}$,
which we will denote as $|{\bf a}|$. For the reader’s convenience, in the following couple of technical lemmas we prove each step in this iterative argument.

\begin{Lemme}[The initial case] With the same hypothesis of Proposition \ref{Prop-Regularity-Classical}, for $|{\bf a}|\leq 2$ and for $1\leq \sigma <+\infty$ we have $\partial^{{\bf a}} \u \in \dot{M}^{2\sigma, p\sigma}(\Rt)$.
\end{Lemme} 
\noindent\pv \\
Let $|{\bf a}|=1$. By equation (\ref{Fixed-point-NS-Classical}) we write
\begin{equation}\label{Fixed-point-NS-Classical-der}
  \partial^{{\bf a}}  \u = -(-\Delta)^{-1} \P\big(\text{div} \partial^{{\bf a}}(\u \otimes \u)\big) + (-\Delta)^{-1}\P  (\text{div}\partial^{{\bf a}}(\mathbb{F})),
\end{equation} 
where we will prove that each term on the right-hand side belongs to the space $\dot{M}^{2\sigma,p\sigma}(\Rt)$. For the first term, since $\u \in \dot{M}^{2,p}(\Rt)\cap L^\infty(\Rt)$ by interpolation inequalities we have $\u \in \dot{M}^{2\sigma,p\sigma}(\Rt)$ for $1\leq \sigma <+\infty$. Then, by H\"older inequalities we obtain $\u \otimes \u \in \dot{M}^{2\sigma,p\sigma}(\Rt)$. As $|{\bf a}|=1$ remark that the operator $(-\Delta)^{-1} \P\big(\text{div} \partial^{{\bf a}}(\cdot)$ writes down as a linear combination of Riesz transforms and by their continuity in the Morrey spaces \cite[Lemme $4.2$]{kato1992p}  we get that $(-\Delta)^{-1} \P\big(\text{div} \partial^{{\bf a}}(\u \otimes \u)\big) \in \dot{M}^{2\sigma,p\sigma}(\Rt)$. Similarly, by (\ref{Hyp-Force}) we get $(-\Delta)^{-1}\P  (\text{div}\partial^{{\bf a}}(\mathbb{F})) \in \dot{M}^{2\sigma,p\sigma}(\Rt)$. 

\medskip

Let $|{\bf a}|=2$. We get back to the expression  (\ref{Fixed-point-NS-Classical-der}) and following  same ideas we can prove that each term on the right-hand side belong to $\dot{M}^{2\sigma,p\sigma}(\Rt)$. We just mention that for the first term we write ${\bf a}={\bf a}_1+{\bf a}_2$, with $|{\bf a}_1|=|{\bf a}_2|=1$. Then, to handle the term ${\partial^{{\bf a_2}}}(\u \otimes \u)$, for $i,j=1,2,3$ we  write $\partial_{i}(u_i u_j)=(\partial_i u_i)u_j+u_i (\partial_i u_j)$ and we use the information $\partial^{{\bf a}_1} \u \in \dot{M}^{2\sigma,p\sigma}(\Rt)$, $\u \in L^{\infty}(\Rt)$ to obtain that $(-\Delta)^{-1} \P\big(\text{div} \partial^{{\bf a}}(\u \otimes \u)\big) \in \dot{M}^{2\sigma,p\sigma}(\Rt)$. 
\finpv

\medskip

\begin{Lemme}[The iterative process] With the same hypothesis of Proposition \ref{Prop-Regularity-Classical}, for $1\leq m \leq k$ and for $|{\bf a}|\leq m$ assume that $\partial^{{\bf a}}\u \in \dot{M}^{2\sigma,p\sigma}(\Rt)$ (with $1\leq \sigma<+\infty$). Then it holds for $|{\bf a}|=k+2$. 
\end{Lemme} 
\noindent\pv\\
Let $|{\bf a }|=k+1$. As before, we must verify that each term on the right-hand side of equation (\ref{Fixed-point-NS-Classical-der}) belongs to $\dot{M}^{2\sigma,p\sigma}(\Rt)$.  To handle the first term, we split ${\bf a}={\bf a}_1+{\bf a}_2$, with $|{\bf a }_1|=1$ and ${|\bf a}_2|=k$, and we  write 
\[ (-\Delta)^{-1} \P\big(\text{div} \partial^{{\bf a}}(\u \otimes \u)\big)= (-\Delta)^{-1} \P\big(\text{div} \partial^{{\bf a}_1} \partial^{{\bf a}_2}(\u \otimes \u)\big). \]
Then, for the term $\partial^{{\bf a}_2}(\u \otimes \u)$, we shall verify that  $\partial^{{\bf a}_2}(\u \otimes \u) \in \dot{M}^{2\sigma, p\sigma}(\Rt)$ for any $\sigma \geq 1$. Indeed, by the classical Leibniz rule (for simplicity we omit the constants) we write
\[ \partial^{{\bf a}_2}(u_i u_j)= \sum_{|{\bf b}| \leq k} \partial^{{\bf b}} u_i \, \partial^{{\bf a}_2 - {\bf b}} u_j=u_i \partial^{{\bf a}_2}u_j + \sum_{1\leq |{\bf b}| \leq k-1} \partial^{{\bf b}} u_i \, \partial^{{\bf a}_2 - {\bf b}} u_j+ \sum_{|{\bf b}|=k} \partial^{\bf b}u_i u_j.\]
Recall that by our hypothesis, for all multi-indice $|{\bf a}|\leq k$ we have $\partial^{\bf a} \u \in \dot{M}^{2\sigma, p\sigma}(\Rt)$. With this information and that fact that  $\u \in L^{\infty}(\Rt)$ we get 
\[ u_i \partial^{{\bf a}_2}u_j, \ \sum_{|{\bf b}|=k} \partial^{\bf b}u_i u_j \in \dot{M}^{2\sigma, p\sigma}(\Rt).\]
Moreover, using again the information $\partial^{\bf a} \u \in \dot{M}^{2\sigma, p\sigma}(\Rt)$, we can write $\partial^{\bf a} \u \in \dot{M}^{2(2\sigma), p(2\sigma)}(\Rt)$ since it holds for any $\sigma \geq 1$. We then obtain
\[ \sum_{1\leq |{\bf b}| \leq k-1} \partial^{{\bf b}} u_i \, \partial^{{\bf a}_2 - {\bf b}} u_j \in  \dot{M}^{2\sigma, p\sigma}(\Rt).\]
Once we have $\partial^{{\bf a}_2}(\u \otimes \u) \in \dot{M}^{2\sigma, p\sigma}(\Rt)$, we conclude that 
$(-\Delta)^{-1} \P\big(\text{div} \partial^{{\bf a}}(\u \otimes \u)\big) \in \dot{M}^{2\sigma, p\sigma}(\Rt)$. For the second term, by our assumption (\ref{Hyp-Force}) we have $(-\Delta)^{-1}\P  (\text{div}\partial^{{\bf a}}(\mathbb{F})) \in \dot{M}^{2\sigma, p \sigma}(\Rt)$. We thus get $\partial^{\bf a} \u \in \dot{M}^{2\sigma, p \sigma}(\Rt)$ for $|{\bf a}|=k+1$

\medskip 

Let $|{\bf a}|=k+2$. Once we have the information above at our disposal, we just repeat again previous arguments to conclude that  $\partial^{\bf a} \u \in \dot{M}^{2\sigma, p \sigma}(\Rt)$ for $|{\bf a}|=k+2$. \finpv 

\medskip

The fact that $\partial^{{\bf a}}P \in \dot{M}^{2\sigma,p\sigma}(\Rt)$ with $|{\bf a}|\leq k+1$ is a direct consequence of identity (\ref{Pressure}) (with $\f=\text{div}(\mathbb{F})$) and our assumption (\ref{Hyp-Force}) on the tensor $\mathbb{F}$. 

\medskip

Finally, by the continuous  embedding $\dot{M}^{2,p}(\Rt)\subset \dot{M}^{1,p}(\Rt)$, and the following result:
\begin{Lemme}[Proposition $3.4$ of \cite{giga1989p}]\label{Lem-U4} Let $f \in \mathcal{S}'(\Rt)$ such that $\vec{\nabla} f \in \dot{M}^{1,p}(\Rt)$, with $p>3$. There exists a constant $C>0$ such that for all  $x,y\in \Rt$ we have $\ds{\vert f(x)-f(y) \vert \leq C\, \Vert \vec{\nabla} f \Vert_{\dot{M}^{1,p}} \, \vert x-y \vert^{1-\frac{3}{p}}}$. 
\end{Lemme}
We obtain that  $\u \in \mathcal{C}^{k+1,s}(\Rt)$ and $P\in \mathcal{C}^{k,s}(\Rt)$ with $s=1-\frac{3}{p}$. Thus, we conclude the proof of Theorem \ref{Th-Reg-Classical}. \finpv 


\nocite{*}
\bibliographystyle{siam}  
\bibliography{references}

\vspace{2cm}

\end{document}